\documentclass[a4paper,12pt]{article}
\usepackage{amssymb}
\usepackage{amsthm}
\usepackage{amsxtra}
\usepackage{amscd}
\usepackage{amsfonts}

\def\ups{\upsilon}

\begin{document}

\def\sect{\section}

\newtheorem{thm}{Theorem}[section]
\newtheorem{cor}[thm]{Corollary}
\newtheorem{lem}[thm]{Lemma}
\newtheorem{prop}[thm]{Proposition}
\newtheorem{propconstr}[thm]{Proposition-Construction}

\theoremstyle{definition}
\newtheorem{para}[thm]{}
\newtheorem{ax}[thm]{Axiom}
\newtheorem{conj}[thm]{Conjecture}
\newtheorem{defn}[thm]{Definition}
\newtheorem{notation}[thm]{Notation}
\newtheorem{rem}[thm]{Remark}
\newtheorem{remark}[thm]{Remark}
\newtheorem{question}[thm]{Question}
\newtheorem{example}[thm]{Example}
\newtheorem{problem}[thm]{Problem}
\newtheorem{excercise}[thm]{Exercise}
\newtheorem{ex}[thm]{Exercise}

\def\Bbb{\mathbb}
\def\cal{\mathcal}
\def\mL{{\mathcal L}}
\def\mC{{\mathcal C}}
\def\mT{{\mathcal T}}
\overfullrule=0pt

\def\tV{\widetilde{V}}
\def\si{\sigma}
\def\prf{\smallskip\noindent{\it        Proof}. }
\def\call{{\cal L}}
\def\nat{{\Bbb  N}}
\def\la{\langle}
\def\ra{\rangle}
\def\inv{^{-1}}
\def\ld{{\rm    ld}}
\def\trdeg{{tr.deg}}
\def\dim{{\rm   dim}}
\def\th{{\rm    Th}}
\def\rest{{\lower       .25     em      \hbox{$\vert$}}}
\def\ch{{\rm    char}}
\def\zee{{\Bbb  Z}}
\def\conc{^\frown}
\def\acl{acl_\si}
\def\cls{cl_\si}
\def\cals{{\cal S}}
\def\mult{{\rm  Mult}}
\def\calv{{\cal V}}
\def\aut{{\rm   Aut}}
\def\ffi{{\Bbb  F}}
\def\ffiti{\tilde{\Bbb          F}}
\def\degs{deg_\si}
\def\calx{{\cal X}}
\def\gal{{\cal G}al}
\def\cl{{\rm cl}}
\def\loc{{\rm locus}}
\def\calg{{\cal G}}
\def\calq{{\cal Q}}
\def\calr{{\cal R}}
\def\caly{{\cal Y}}
\def\aff{{\Bbb A}}
\def\cali{{\cal I}}
\def\calu{{\cal U}}
\def\epsilon{\varepsilon} 
\def\Uu{{\cal U}}
\def\rat{{\Bbb Q}}
\def\ga{{\Bbb G}_a}
\def\gm{{\Bbb G}_m}
\def\cee{{\Bbb C}}
\def\ree{{\Bbb R}}
\def\frob{{\rm Frob}}
\def\Frob{{\rm Frob}}
\def\fix{{\rm Fix}}
\def\Uu{{\cal U}}
\def\proj{{\Bbb P}}
\def\sym{{\rm Sym}}
 
\def\dcl{{\rm dcl}}
\def\calm{{\mathcal M}}

\font\helpp=cmsy5
\def\semdp
{\hbox{$\times\kern-.23em\lower-.1em\hbox{\helpp\char'152}$}\,}

\def\dnfo{\,\raise.2em\hbox{$\,\mathrel|\kern-.9em\lower.35em\hbox{$\smile$}
$}}
\def\dnf#1{\lower1em\hbox{$\buildrel\dnfo\over{\scriptstyle #1}$}}
\def\dfo{\;\raise.2em\hbox{$\mathrel|\kern-.9em\lower.35em\hbox{$\smile$}
\kern-.7em\hbox{\char'57}$}\;}
\def\df#1{\lower1em\hbox{$\buildrel\dfo\over{\scriptstyle #1}$}}        
\def\stab{{\rm Stab}}
\def\qfcb{\hbox{qf-Cb}}

%%%%%%Udi's macros
\newcommand{\nc}{\newcommand}
\nc{\renc}{\renewcommand}
\nc{\ssec}{\subsection}
\nc{\sssec}{\subsubsection}
\nc{\on}{\operatorname}

\nc\ol{\overline}
\nc\wt{\widetilde}
\nc\wh{\widehat}
\nc\tboxtimes{\wt{\boxtimes}}

\emergencystretch=2cm

\nc{\Aa}{{\mathbb{A}}}
 \nc{\Gg}{{\mathbb{G}}}  \def\gg{\mathbb{g}}
\nc{\Hh}{{\mathbb{H}}}
 \nc{\Nn}{{\mathbb{N}}}
\nc{\Pp}{{\mathbb{P}}}
\nc{\Rr}{{\mathbb{R}}}
\nc{\BV}{{\mathbb{V}}}
\nc{\BW}{{\mathbb{W}}}
\nc{\Zz}{{\mathbb{Z}}}
\nc{\Qq}{{\mathbb{Q}}}
\nc{\Ss}{{\mathbb{S}}}
\nc{\Cc}{{\mathbb{C}}}

\nc{\CA}{{\mathcal{A}}}
\nc{\CB}{{\mathcal{B}}}

\nc{\CE}{{\mathcal{E}}}
\nc{\CF}{{\mathcal{F}}}
\nc{\CG}{{\mathcal{G}}}
\nc{\CL}{{\mathcal{L}}}
\nc{\CC}{{\mathcal{C}}}
\nc{\CM}{{\mathcal{M}}}
\def\Mm{\CM}
\nc{\CN}{{\mathcal{N}}}
\nc{\Oo}{{\mathcal{O}}}
\nc{\CP}{{\mathcal{P}}}
\nc{\CQ}{{\mathcal{Q}}}
\nc{\CR}{{\mathcal{R}}}
\nc{\CS}{{\mathcal{S}}}
\nc{\CT}{{\mathcal{T}}}
\nc{\CU}{{\mathcal{U}}}
\nc{\CV}{{\mathcal{V}}}
\nc{\CK}{{\mathcal{K}}}
\nc{\CW}{{\mathcal{W}}}
\nc{\CZ}{{\mathcal{Z}}}

\nc{\cM}{{\check{\mathcal M}}{}}
\nc{\csM}{{\check{\mathcal A}}{}}
\nc{\oM}{{\overset{\circ}{\mathcal M}}{}}
\nc{\obM}{{\overset{\circ}{\mathbf M}}{}}
\nc{\oCA}{{\overset{\circ}{\mathcal A}}{}}
\nc{\obA}{{\overset{\circ}{\mathbf A}}{}}
\nc{\ooM}{{\overset{\circ}{M}}{}}
\nc{\osM}{{\overset{\circ}{\mathsf M}}{}}
\nc{\vM}{{\overset{\bullet}{\mathcal M}}{}}
\nc{\nM}{{\underset{\bullet}{\mathcal M}}{}}
\nc{\oD}{{\overset{\circ}{\mathcal D}}{}}
\nc{\obD}{{\overset{\circ}{\mathbf D}}{}}
\nc{\oA}{{\overset{\circ}{\mathbb A}}{}}
\nc{\op}{{\overset{\bullet}{\mathbf p}}{}}
\nc{\cp}{{\overset{\circ}{\mathbf p}}{}}
\nc{\oU}{{\overset{\bullet}{\mathcal U}}{}}
\nc{\oZ}{{\overset{\circ}{\mathcal Z}}{}}
\nc{\ofZ}{{\overset{\circ}{\mathfrak Z}}{}}
\nc{\oF}{{\overset{\circ}{\fF}}}

\nc{\fa}{{\mathfrak{a}}}
\nc{\fb}{{\mathfrak{b}}}
\nc{\fg}{{\mathfrak{g}}}
\nc{\fgl}{{\mathfrak{gl}}}
\nc{\fh}{{\mathfrak{h}}}
\nc{\fj}{{\mathfrak{j}}}
\nc{\fm}{{\mathfrak{m}}}
\nc{\fn}{{\mathfrak{n}}}
\nc{\fu}{{\mathfrak{u}}}
\nc{\fp}{{\mathfrak{p}}}
\nc{\fr}{{\mathfrak{r}}}
\nc{\fs}{{\mathfrak{s}}}
\nc{\fsl}{{\mathfrak{sl}}}
\nc{\hsl}{{\widehat{\mathfrak{sl}}}}
\nc{\hgl}{{\widehat{\mathfrak{gl}}}}
\nc{\hg}{{\widehat{\mathfrak{g}}}}
\nc{\chg}{{\widehat{\mathfrak{g}}}{}^\vee}
\nc{\hn}{{\widehat{\mathfrak{n}}}}
\nc{\chn}{{\widehat{\mathfrak{n}}}{}^\vee}

\nc{\fA}{{\mathfrak{A}}}
\nc{\fB}{{\mathfrak{B}}}
\nc{\fD}{{\mathfrak{D}}}
\nc{\fE}{{\mathfrak{E}}}
\nc{\fF}{{\mathfrak{F}}}
\nc{\fG}{{\mathfrak{G}}}
\nc{\fK}{{\mathfrak{K}}}
\nc{\fL}{{\mathfrak{L}}}
\nc{\fM}{{\mathfrak{M}}}
\nc{\fN}{{\mathfrak{N}}}
\nc{\fP}{{\mathfrak{P}}}
\nc{\fU}{{\mathfrak{U}}}
\nc{\fV}{{\mathfrak{V}}}
\nc{\fZ}{{\mathfrak{Z}}}

\nc{\bb}{{\mathbf{b}}}
\nc{\bc}{{\mathbf{c}}}
\nc{\bd}{{\mathbf{d}}}
\nc{\be}{{\mathbf{e}}}
\nc{\bj}{{\mathbf{j}}}
\nc{\bn}{{\mathbf{n}}}
\nc{\bp}{{\mathbf{p}}}
\nc{\bq}{{\mathbf{q}}}
\nc{\bF}{{\mathbf{F}}}
\nc{\bu}{{\mathbf{u}}}
\nc{\bv}{{\mathbf{v}}}
\nc{\bx}{{\mathbf{x}}}
\nc{\bs}{{\mathbf{s}}}
\nc{\by}{{\mathbf{y}}}
\nc{\bw}{{\mathbf{w}}}
\nc{\bA}{{\mathbf{A}}}
\nc{\bK}{{\mathbf{K}}}
\nc{\bI}{{\mathbf{I}}}
\nc{\bB}{{\mathbf{B}}}
\nc{\bG}{{\mathbf{G}}}
\nc{\bC}{{\mathbf{C}}}
\nc{\bD}{{\mathbf{D}}}
\nc{\bP}{{\mathbf{P}}}
\nc{\bH}{{\mathbf{H}}}
\nc{\bM}{{\mathbf{M}}}
\nc{\bN}{{\mathbf{N}}}
\nc{\bV}{{\mathbf{V}}}
\nc{\bU}{{\mathbf{U}}}
\nc{\bR}{{\mathbf{R}}}
\nc{\bL}{{\mathbf{L}}}
\nc{\bT}{{\mathbf{T}}}
\nc{\bW}{{\mathbf{W}}}
\nc{\bX}{{\mathbf{X}}}
\nc{\bY}{{\mathbf{Y}}}
\nc{\bZ}{{\mathbf{Z}}}
\nc{\bS}{{\mathbf{S}}}
\nc{\bQ}{{\mathbf{Q}}}

\nc{\sA}{{\mathsf{A}}}
\nc{\sB}{{\mathsf{B}}}
\nc{\sC}{{\mathsf{C}}}
\nc{\sD}{{\mathsf{D}}}
\nc{\sF}{{\mathsf{F}}}
\nc{\sG}{{\mathsf{G}}}
\nc{\sK}{{\mathsf{K}}}
\nc{\sM}{{\mathsf{M}}}
\nc{\sO}{{\mathsf{O}}}
\nc{\sQ}{{\mathsf{Q}}}
\nc{\sP}{{\mathsf{P}}}
\nc{\sZ}{{\mathsf{Z}}}
\nc{\sfp}{{\mathsf{p}}}
\nc{\sr}{{\mathsf{r}}}
\nc{\sg}{{\mathsf{g}}}
\nc{\sff}{{\mathsf{f}}}
\nc{\sfb}{{\mathsf{b}}}
\nc{\sfc}{{\mathsf{c}}}
\nc{\sd}{{\ltimes}}

\nc{\tA}{{\widetilde{\mathbf{A}}}}
\nc{\tB}{{\widetilde{\mathcal{B}}}}
\nc{\tg}{{\widetilde{\mathfrak{g}}}}
\nc{\tG}{{\widetilde{G}}}
\nc{\TM}{{\widetilde{\mathbb{M}}}{}}
\nc{\tO}{{\widetilde{\mathsf{O}}}{}}
\nc{\tU}{{\widetilde{\mathfrak{U}}}{}}
\nc{\TZ}{{\tilde{Z}}}
\nc{\tx}{{\tilde{x}}}
\nc{\tq}{{\tilde{q}}}

\nc{\tfP}{{\widetilde{\mathfrak{P}}}{}}
\nc{\tz}{{\tilde{\zeta}}}
\nc{\tmu}{{\tilde{\mu}}}

  \nc{\Ob}{{\mathop{\operatorname{\rm Ob}}}}
  \nc{\Sym}{{\mathop{\operatorname{\rm Sym}}}}
   \nc{\Aut}{{\mathop{\operatorname{\rm Aut}}}}
 \nc{\Spec}{{\mathop{\operatorname{\rm Spec}}}}
  \nc{\spec}{{\mathop{\operatorname{\rm Spec}}}}
\nc{\Ker}{{\mathop{\operatorname{\rm Ker}}}}
 \nc{\dom}{{\mathop{\operatorname{\rm dom}}}}
\nc{\End}{{\mathop{\operatorname{\rm End}}}}
 \nc{\Hom}{\on{\Hom}}
 \nc{\GL}{{\mathop{\operatorname{\rm GL}}}}
 \nc{\Id}{{\mathop{\operatorname{\rm Id}}}}
 \nc{\rk}{{\mathop{\operatorname{\rm rk}}}} 
 \nc{\length}{{\mathop{\operatorname{\rm length}}}}
\nc{\supp}{{\mathop{\operatorname{\rm supp}}}}
\nc{\val}{{\rm val}}
%\nc{\res}{{\rm res}}
\nc{\res}{{\mathop{\operatorname{\rm res}}}}
\def\ind#1#2{ {#1} {\downarrow} {#2} }  
\def\Ind#1#2#3{{#1} {\downarrow}_{#3} {#2} }  
\def\domn#1#2#3{ (#1 \leq_{dom_{#3}} #2 ) }
\def\tensor{{\otimes}}
\def\meet{\cap}
\def\union{\cup}
\def\Union{\bigcup}
\def\si{\sigma}
\def\g{\gamma}
\def\G{\Gamma}
\def\Sum{\Sigma}
\def\<{\begin}
 \def\>{\end}
\def\m{\setminus}

 \def\AFD{(FD)\,}
\def\Aom{(FD$_{\omega}$)\,}
\nc{\seq}[1]{\stackrel{#1}{\sim}}
\def\dd#1{\frac{\partial}{\partial X_{#1}}}   \def\ddt{\frac{\partial}{\partial t}}
\def\invv#1{{{#1}^{-1}}}
\def\inv{^{-1}}
\def\claim#1{{\noindent \bf Claim #1.\ }}
\def\Claim{{\noindent \bf Claim.\ }}
\def\beq#1{\begin{equation} \label{#1}  }
\def\eeq{\end{equation}}
\def\normal{\trianglelefteq}
\def\Uu{\mathcal U}
\def\iso{\simeq}

\def\bt{\beta}
\def\pv{\hfill $\Box$}
\def\eprf{\hfill $\Box$}
 
\def\acl{\mathop{\rm acl}\nolimits}
 \def\dcl{\mathop{\rm dcl}\nolimits}
\def\OK{{\mathcal O}_K}
\def\liminv{\underset{\longleftarrow}{lim}\,} 
\def\liminvi{\underset{\longleftarrow i}{lim}\,} 

\def\lbl#1{    \label{#1}  }
\def\a{\alpha}
\def\revo{\backslash}
\def\ba{\bar{a}}
\def\bba{{\mathbf a}}
\def\k{{\rm k}}
 
\def\tlt{\times \ldots \times}

  \def\Rr{\ree}

\def\frac#1#2{{#1\over #2}}

\def\lam{\lambda}
\def\e{\epsilon}

\def\vlabel{\label}

\def\udirem{\smallskip\noindent {\bf Remark}. }

\title{Difference fields and descent in algebraic dynamics - I  \thanks{Partially supported by project ANR-06-BLAN-0183, by Modnet: MRTN-CT-2004-512234,  and by Israel Science Foundation grant 1048-07}}

\author{Zo\'e Chatzidakis \ { (CNRS - Universit\'e Paris 7)}
\and
Ehud Hrushovski   
(The Hebrew University at Jerusalem)    
}
%\centerline{\today}
\maketitle
\def\hK{{K}}
\def\hL{{\widehat{L}}}

\abstract{

We draw a  connection between the model-theoretic notions of modularity (or one-basedness), orthogonality and internality, as applied to difference fields, and questions of descent in  in algebraic dynamics.   In particular we prove in any dimension a strong dynamical version of Northcott's theorem for function fields, answering a question of Szpiro and Tucker and generalizing a theorem of Baker's for the projective line.  

   The  paper comes in three parts.  This first  part contains an exposition some of the main results of the model theory of difference fields, and their immediate connection to questions of descent in algebraic dynamics.   We present the model-theoretic notion of internality in a context that does not require a universal domain with quantifier-elimination.   We also note a version of canonical heights that applies well beyond polarized algebraic dynamics.
Part II sharpens the structure theory to arbitrary base fields and constructible maps  where in part I we emphasize finite base change and correspondences.   Part III will 
include precise structure  theorems related to the Galois theory considered here, and will 
enable a sharpening of the descent results for non-modular dynamics. }
 
\section{Introduction} 

Algebraic dynamics studies algebraic varieties $V$   with an endomorphism $\phi$.  
Most questions reduce to the case that $\phi$ is dominant, and we will assume
this.  The base field is taken to be a number field or a function field $F$; given an extension field $K$ of $F$,
one has the set $V(K)$, and a self-map $\phi: V(K) \to V(K)$.  
 
The same objects $\bV=(V,\phi)$ arise in model theory in a somewhat different way.  Here we take not only $K$ but also an endomorphism $\si$ of $K$, and associate to $\bV$ the 
fixed points of $\phi$ twisted by $\si$, i.e. the set $\bV(K,\si) = \{x \in V(K): \si(x) = \phi(x)\}$.  The pair
$(K,\si)$ is called a difference field; the functor $(K,\si) \mapsto \bV(K,\si)$ suffices to 
recover $\bV$.  This point of view makes just as much sense when one takes
$\phi$ to be any correspondence on $V$, and the theory is carried out in that generality;
we will restrict attention in the introduction to the case that $\phi$ is a rational function. 

   Beyond the coincidence of the objects, the two subjects  diverge.  Much of the 
richness of algebraic dynamics comes from the arithmetic of the field,
and in particular 
  the interaction of the dynamics with valuations and absolute
values on $K$.    But the  model theory of valued difference fields   is still in its infancy, 
and the model theory of global fields is not yet conceived.   Current
model-theoretic results thus concern purely geometric aspects of
algebraic dynamics.  On the other hand, in a number of ways that will be
detailed below, the treatment is more general; and in particular 
there is an effort to find precise dividing lines among quite
differently behaved dynamics, in any dimension.  It seems possible that
the general properties that arise from this study can  be useful 
in algebraic dynamics, and the purpose of this paper
is to bring them out.   

Our motivation lay in a question of   Szpiro and Tucker concerning descent for algebraic dynamics, arising out of  Northcott's theorem for dynamics over function fields.  
The
question was originally formulated for polarized dynamical systems in terms of canonical heights; Szpiro translated it into the language of limited sets, where it admits a much more general form, and it is in this form that we solve the problem. 
A subset of the function field that can be parametrized by a constructible set over the base
field will be called {\em limited}; cf. \S \ref{limited}.\footnote{ ``bounded'' may be a more
common term;  \cite{langneron} speak of points ``belonging to a finite number of algebraic families".}   
For algebraic dynamics $\phi$ on $\Pp^1$ with $\deg(\phi)>1$, Baker showed
that no infinite orbit can be contained in a limited set, unless the dynamics is isotrivial
i.e. comes from a dynamics defined over the base field.  Generalizing this to higher
dimensions  requires first of all a definition:  what is the right analogue, for higher 
  dimensions, of the hypothesis that $\deg(\phi)>1$?   Of isotriviality?  For 
 Abelian varieties, the notion of a simple component is clear, and
 isotriviality needs to be refined to a notion of a   {\em trace} of a
 family of Abelian varieties, roughly speaking 
 the sum of simple components that stay constant in the family.  What is the analogue when one
 generalizes to algebraic dynamics?  The model theoretic description of the category
 of algebraic dynamics gives clear answers; using them we   sharpen and generalize 
 Baker's theorem to arbitrary dimension.  As Daniel Bertrand pointed out to us, in
the case of Abelian varieties with a multiplication dynamics, our result is classical and due
  to Lang and N\'eron; cf. \cite{lang} Chapter 6, Theorem 5.4.
 %:add something somewhere about descent of Abelian varieties.

  Model theory usually begins by amalgamating all relevant structures into a
homogeneous one, a universal domain $\Uu = (\Uu,\si)$.   Then $\bV(\Uu)$ contains all the information in the functor $(K,\si) \mapsto \bV(K,\si)$.  Now difference fields
admit amalgamation over algebraically closed subfields, but not over arbitrary
subfields.  This obliges us   to consider more than one universal domain; $\Uu$
will be determined by its prime field $F$, and   the   conjugacy class of $\si$ in
  $Aut(F^{alg}/F)$.     Moreover, the family of  sets $\bV(\Uu)$ (or their Boolean combinations)
  is not closed under projections.  However,    the fact that amalgamation holds over algebraicallly closed subfields means
  that it is not necessary to look further than projections under finite maps.   By
a {\em definable subset} of $V(\Uu)$ we mean one defined by some formula
formed out of diference equations and inequations by means of the logical operators
$(\exists x),(\forall y), \wedge, \neg$. The following Proposition  was proved for $V(F)$,
the set of points of $V$ in a  pseudo-finite field $F$, by Ax.  Van den Dries understood that Ax's theorem can be interpreted
as applying to $\bV(\Uu)$ where $V$ is a variety with the trivial dynamics, suggesting
that the same form of quantifier elimination may be valid more generally.    This was carried out  in \cite{mac},  as well as in \cite{CH}.   Here we state it for algebraic
dynamics, though a version valid for any difference variety is available

\<{prop}  \vlabel{qe} Every definable subset of $\bV(\Uu)$ is a finite Boolean combination
of sets  of the form $f(\bW(\Uu))$, where $\bW=(W,\phi') \in AD_K$ and $f: W \to V$
is a quasi-finite morphism of varieties, with $\phi \circ f=f \circ \phi'$. \>{prop} 
  
We call $\bW$ a {\em finite cover} of $\bV$  (or rather of the closure of the image of $f$.)    Such finite covers cannot be avoided;
indeed if $\deg(f)>1$, possibly after passing from $\si$ to $\si^n$, the set $f(\bW(\Uu))$
is never a Boolean combination of quantifier-free definable sets; moreover if
$f(\bW(\Uu))$ is a Boolean combination of sets $f_i(\bW_i(\Uu))$, then the finite cover
$W \to V$ is itself a quotient of the fiber product of the covers $W_i$.

 The theory of $\Uu$ is {\em simple}, and stable on a quantifier-free level.   The methods
and definitions of stability theory apply to $\bV(\Uu)$:  stability, modularity, internality to fixed fields, analysis in terms of minimal types.   We would like to view these as properties of $\bV$;
but the embedding into $\Uu$
is not canonically determined by the geometric data $\bV$, and moreover, even if the variety
$V$ is absolutely irreducible, the embedding 
  splits the generic point of $\bV$ into a possibly infinite number of generic types.    
%
 %induces additional structure on $\bV$; namely a choice of a conjugacy class of a lifting of the endomorphism induced
%by $\phi$ on the function field of $V$ to the algebraic closure of that function field.  On the other hand, the morphisms allowed model-theoretically are more general; we
%allow multi-valued correspondences and not only morphisms.  
We show here that the above 
do not depend on the embedding, or the choice of generic type, but are really geometric properties
of $\bV$. Moreover, the 
 decomposition theory can be carried out using rational maps.  
 In particular the notion of descent depends little on whether one chooses rational, constructible or  multi-valued morphisms.  
 We do this in two ways:  in \S 2, we develop the basic theory of internality from
scratch for the category of difference fields, without a preliminary amalgamation into 
a universal domain.   In this way the theory is more general, applying for instance to pairs
of commuting automorphisms.  
In Part II of this paper, we use the usual model-theoretic language
but show a posteriori the independence of the embedding; in this way we do not lose
the  intuitions associated with the model-theoretic viewpoint.   In both approaches, a
  certain weakening of amalgamation plays  an   essential role.  This weakening
  is associated to what one might loosely call
``Shelah's reflection principle".   Roughly speaking, stable interactions between a type $P$ 
and external elements   can be read off from inside $P$ itself.  Algebraically, the difficulty
with difference algebra is that difference fields do not admit amalgamation; however
if $K \leq L$ is an extension of difference fields, there is always a canonical amalgam of $L$
with {\em itself} over $K$, identifying the algebraic part; using a version of the reflection principle, we show that this suffices for a definable Galois theory.

In this introduction, we will restrict attention to  the category of algebraic dynamics.   The results are valid in greater generality, for difference varieties that arise from correspondences rather than rational maps.   
 
We are grateful to the referee for very helpful comments.

\para{\noindent \bf Modularity.}  
Let $K$ be a field, and $AD_K$ the category of pairs $\bV=(V,\phi)$ where $V$ is an irreducible variety over $K$, and $\phi: V \to V$ a  dominant rational map.   In this paper, we will only consider ``birational'' or ``generic''  questions, i.e. we allow ourselves to  ignore  any given  lower dimensional subvariety  of $V$.  We correspondingly    take a morphism
$(V,\phi) \to (V',\phi')$ to be a dominant rational map $f: V \to V'$, such that
$\phi' \circ f = f \circ \phi$.   
  In particular   $(V,\phi)$ is considered  isomorphic  to 
$(V',\phi')$ if there exists an isomorphism $V \setminus U \to V' \setminus U'$ for some 
lower-dimensional $U,U'$, commuting (whenever defined) with $\phi$, $\phi'$.   At this level we could dispense
with varieties altogether, and speak of their function fields instead.  However both the proofs
of the statements, and intended future developments, require the geometric viewpoint.

If $\bU=(U,\phi) \in AD_K$,
we let $AD_{\bU}$ be the category of $AD_K$-morphisms $\bV \to \bU$.   
Note that a fiber of $V \to U$ is a subvariety of $V$, which is not in
general invariant under the dynamics. 
%so a fiber of $\bV \to \bU$ is not an algebraic dynamics in the above sense.  
However if we generalize the notion of algebraic dynamics to difference
fields, we can speak of the fibers of $\bV \to \bU$.  In general, given
a difference 
field $\hK = (K,\si)$, let $AD_\hK$ be the  category of triples
$(V,V^\si,\psi)$, with $\psi$ a dominant rational map $V \to V^\si$.  Here $V^\si$ is obtained
from $V$ by applying $\si$ to the coefficients.
  In
our case,  
 $\phi$ induces an endomorphism $\si$ of the function field  of $U$,
 making it into a difference field $\hK$.  The generic fiber of $\bV \to \bU$
 can be understood as an object of $AD_\hK$, and indeed $AD_\bU$ is   isomorphic to $AD_{\hK}$.   In this introduction we will stay with the geometric language.
 
%AD_K
For $\bV \in AD_K$, the irreducible components of $\bV \times \bU$, with
the projection 
maps to $\bU$, are elements of $AD_\bU$.   If  $\bV \times \bU$ is irreducible, we denote
it $\bV_\bU$.     More generally, if $V$ remains irreducible over a difference field
extension $L$ of $K$, we write $\bV_L$ for the change of basis.

 If $\bW,\bV$ are objects of $AD_K$, we write $U  \leq \bV \times \bW$ 
  if $U$ is an irreducible subvariety
of $V \times W$  such that the projections $U \to V, U \to W$ are $AD_K$-morphisms
with respect to some (unique) $\phi_\bU$; in this case we write also 
$\bU \leq \bV \times \bW$.
 We will sometimes think of $U$ as {\em a
family of difference subvarieties of $\bV$}, indexed by $a \in W$.  Assume that for generic $a \in W$, $U_a$ is absolutely
irreducible  of dimension $l$,  and if $b \neq a$ then $U_a \neq U_b$; then we say that it is 
a family of dimension $\dim(W)$ of $l$-dimensional difference subvarieties of $\bV$.
By a {\em difference subvariety} of $\bV$ we mean a generic fiber of some such   family.

In algebraic geometry there exist $n$-dimensional families of irreducible plane curves
for arbitrarily large $n$.   It is a  fundamental attribute of algebraic dynamics that -with rare exceptions -   the dimension of families is bounded.  To state this precisely we need
to demarcate off the exceptional sub-category,  that of   field-internal dynamics, that behaves like algebraic geometry.

Say $\bV=(V,\phi)$  has constant dynamics if $\phi=Id_V$; {\em periodic}, if for some $n$ we have $\phi^n=Id$;  {\em twisted-periodic} if   $\phi^n$ is a   Frobenius morphism on $V$.\footnote{This requires   $V$ to be defined over a finite field; or in the more general situation
of $AD_{(L,\si)}$ considered below, that $V$ descend to $Fix(\si^{n} Frob^{m})$ for some $n \in \Nn, m \in \Zz$.}

 An object of $AD_\bU$ is called 
 constant (periodic, twisted-periodic)   over $\bU$  if it is isomorphic in $AD_\bU$ to some 
 $\bV  \leq \bQ \times \bU$, with $\bQ$ enjoying the corresponding property.

Call $\bV \in AD_\bU$ 
  {\em field-internal} if  for some $\bW \in AD_\bU$, some   generic component 
  $\bV'$ of $\bV \times_\bU \bW$ is   twisted-periodic in $AD_\bW$.\footnote{We quote a sentence from the referee report regarding the terminology:
the ``notions of ``field-internal," ``fixed-field-internal," and ``field-free" are applied to
algebraic dynamics but they make essential reference to the difference field point of view   \ldots
these  terms were selected to facilitate the difference field interpretation."}
    The untwisted analogue
  will be referred to as {\em fixed-field-internal.}   % \footnote{If $K$ is not algebraically closed then $U \times V$ need not itself be irreducible.}
  
The central example is given by {\em translation varieties}.  Let $G \times X \to X$ be an algebraic group action defined over   $K$.  Let $g: U \to G$.  Let $V= U \times X$,
and define a dynamics by:  $(u,x) \mapsto (\phi_\bU(u), g(u) \cdot x)$.   In this case 
$\bV$ is field-internal over $\bU$.   When  $G=GL_n$ and $X$ is the natural module, this is the object of Van der Put - Singer Picard-Vessiot theory.

   Finally, $\bV $ is {\em field-free}  over $\bU$ if
   it has no field-internal components, i.e if 
  whenever $\bV \to \bW \to \bU'  \to \bU$  factorizes the given morphism $\bV \to \bU$ and $\bW$ is field-internal over $\bU'$, then   $\bW \to \bU$
has finite fibers (generically over $\bU$).     Similarly for fixed-field-free.  For $\bU$ a point, we say $\bV$ is field-free.

\<{prop} \vlabel{mod1} 
\begin{itemize}
\item[(1)] {   If $\bV \leq \bV' \times \bV''$ and $\bV',\bV''$ are
field-free, so is $\bV$.   }
\item[(2)]  {Let $\bV$ be field-free, and let $\bU$ be field-internal.   Let
$R \leq V \times U$ be an irreducible difference subvariety, projecting dominantly to $U,V$. 
 Then $R$ is a   component  of $V \times U$.} \>{itemize}  \>{prop}

(1) is just a closure property of the class of field-free dynamics.  (2) is the quantifier-free part of   {\em orthogonality to fixed fields}:  there is no quantifier-free interaction between field-free and field-internal dynamics.   We must speak of a component since $V \times U$ may not be irreducible, but if it is, the conclusion is that the only possible relation $R$ is the trivial one.  
It suffices here to assume that $\bV$ has no nontrivial field-internal quotients.  
We sketch a model-theoretic proof of a variant of this, 
 illustrating how the solution set $\bV(\Uu)$, and especially the  model-theoretic notion of    induced structure on $\bV(\Uu)$,   show up here.

(2)  reduces to the case that $\bU$ is twisted-periodic, so that each point of
$\bU(\Uu)$ is contained in a difference field generated over $K$ by $Fix(\tau)$,
$\tau=\si^m Frob^l$.  Take $\tau=\si$ for simplicity.   Then  (2) 
 is actually true under the weaker assumption (*):  $\bV$ has no positive-dimensional quotients with constant dynamics.    We sketch the proof.

By (*), for a  generic point $b$ of   $\bV(\Uu)$, we have $K(b) \meet   Fix(\si)
\subseteq K^{alg}$.
%  In fact $K(b) \meet Fix(\si^n) \subseteq K^{alg}$ for any $n$,
%since if $(V_1,\phi_1)$ is a quotient of $\bV$ with periodic dynamics $\phi_1^m=Id$,
%then the set of orbits of $\phi_1$ on $V_1$ is again a quotient of $\bV$, with constant dynamics.
It follows that  $K(b)^{alg} \meet   Fix(\si) \subseteq K^{alg}$ (see Lemma \ref{2-1}.)  
 
On the other hand using quantifier elimination to the level of images of finite maps (cf. Proposition \ref{qe}) one sees that definable
closure is contained in the field-theoretic algebraic closure.   It follows
that every $K$- definable map $\bV(\Uu) \to Fix(\si)$ is generically constant. 
Using   elimination of imaginaries   and stable embeddedness of $Fix(\si)$, it follows that if $c$ is a generic point of $\bV(\Uu)$,
then every $K(c)$-definable subset of $Fix(\si)^n$ is $K^{alg}$-definable:  the code
for such a set being a definable function of $c$.  It follows that every definable
relation on $\bV(\Uu) \times Fix(\si)^n$ is a finite Boolean combination of rectangles $X \times Y$, with $X \subseteq \bV(\Uu)$ and $Y \subseteq Fix(\si)^n$, and of relations whose projection
to $\bV(\Uu)$ is not generic.  From this (2) is immediate.

  In fact   (2) does not exhaust the strength of the model-theoretic orthogonality just proved;  it is stated only for quantifier-free definable sets, whereas 
   we showed  that    there is no interaction using  subsets
of $\bV(\Uu) \times \bU(\Uu)$ that are defined using quantifiers, either:  any
 $\bU(\Uu)$-parametrized family of definable
subsets of $\bV(\Uu)$ must be finite.  Using the remark following Proposition \ref{qe},
we can translate this to a statement about finite covers.  If ${\mathcal F}$ is
a $\bU(\Uu)$-parameterized family of finite covers $f: \bW \to \bV$, then the family
of sets $f(\bW(\Uu))$ must be finite.  From the converse to \ref{qe} it follows
that the elements of ${\mathcal F} $ themselves arise from a finite family of finite covers.
More precisely, we have:

\begin{itemize}
\item[($2^+)$] {   let $L=K(\bU)$,
and let $\bV''$ be  a finite cover of $\bV_{L}$.  Then there exists a finite cover of $\bV''$
of the form $\bV' _{L}$, where $\bV'$ is a 
finite cover  of $\bV$.  } \end{itemize}
 
 This is  Lemma 4.2 of \cite{CHS}.  The same proof applies in the twisted-periodic case; compare 1.11(3) of Part II.

So far we only noted closure and orthogonality properties.  Here is the essential statement:  
 
\<{thm} \vlabel{mod}  Assume $(V,\phi)$ is field-free.   Let $R \leq \bV \times \bQ$ be an irreducible $k$-dimensional family of
$l$-dimensional  irreducible 
  difference subvarieties of $\bV$.  Then $k+l \leq \dim(V)$.    
 
  \>{thm}

 This property is characteristic of theories of modules.  We say $\bV$ is {\em modular} if (3) holds for all  powers $\bV^n$ of $\bV$.   Thus for algebraic dynamics, modularity is equivalent to being field-free.

Modularity is a fundamental dividing line in model theory, and has many equivalent
formulations.  Here we will just note an algebraic equivalent that will be used in Part II.
Call $\bV$ {\em one-based} if for any difference field $L$ extending $K$, any tuple $a$ from  $\bV(L)$, and any tuple $b$ from $L$, the fields $K(a)^{alg},K(b)^{alg}$   are free over
their intersection.  
%Modularity of $\bV$ is equivalent to {\em one-basedness}, which can be stated geometrically
%as follows:   Let $\bU \in AD_K$, and let $R \leq \bU \times \bV$.
%Then there exist morphisms $\bU' \to \bU,  \bV' \to \bV$   with finite fibers, and morphisms
%$\bV' \to \bW$, $ \bU' \to \bW$, such that $R$ is just the image of $U' \times _W V' $ in $U \times V$.   
% It is easy to see that one-basedness implies modularity.  To see the converse,  assume $\bV$ is modular;  view $R$ as a family of difference subvarieties of $\bV$;  passing to   finite cover $\bU'$ of
%$\bU$, we can take them to be  generically irreducible, of dimension $e$;   modularity implies
%that $R$ is a pullback of some $\bV' \leq \bW \times \bV$, with  $\bW $ a quotient of $\bU'$ of dimension $\dim(V)-e$; the projection $\bV' \to \bV$ thus has finite fibers.   
  To see  that modularity implies  one-basedness and modularity,  say $\bV,\bU,\bR \in AD_K$
	have function fields isomorphic to $K(a)_\si,K(b)_\si,K(a,b)_\si$
	respectively.   
	 Let $R_b = \{v \in V: (v,b) \in R \}$, and let $U$ be the irreducible component of $R_b$
	 containing $a$; let  of $K(b')$ be the field of definition of $U$.  Then $b' \in K(b)^{alg}$.
	By the modular rank inequality we have $tr.deg._K( K(b'))+ tr. deg._{K(b')} K(a,b') \leq
	tr. deg._K K(a)$; so $b' \in K(a)^{alg}$.  Thus $K(a),K(b)$ are free over
	the intersection of their algebraic closures.  The converse is proved similarly.   We will reserve
	 the term {\em one-based} to the solution set $\bV(\Uu)$, and {\em modular} to the 
	 geometric data $(V,\phi)$ itself.

 Theorem \ref{mod}    expresses a dichotomy between modular and field-like behavior in algebraic dynamics.
This strengthens \cite{CHP}, where in effect modularity is proved assuming every finite cover
of $(V,\phi)$ is field-free.     

 The dichotomy is an expression of a general   philosophy of Zilber's, and indeed
 it is through this general principle that the theorem is proved in \cite{CHP}.   The proof is 
 easy to explain in the case of algebraic dynamics over a field.  We say
  $\bV=(V,\phi)$ is {\em primitive}  if there is no   $f: \bV \to \bU$, $0< \dim(U) < 
\dim(V)$.   
 To avoid technicalities,   assume $\phi$ is an endomorphism of a smooth variety $V$, and that for any
 $m \geq 1$  there are no $\phi^m$-invariant proper, infinite subvarieties, 
 nor any dominant, equivariant $f: (V,\phi^m) \to (U,\psi)$ with  $0< \dim(U) < 
\dim(V)$.  
Consider $V$ as a structure, with $n$-ary relations given by   subvarieties of $V^n$ 
 left invariant by a power of $\phi$.  Then it is very easy to check that we have a Zariski geometry, i.e. a structure with a topology satisfying the basic properties of Zariski closed
 sets.    \footnote{
 The main point is that if $U_1,U_2$ are periodic then so is each component of $U_1 \meet U_2$.    It is in this that the notion of ``periodic'' behaves better than "invariant".
 The Boolean combinations of periodic subvarieties are also closed under projections, since
 both the closure and the boundary of the  projection are invariant for the same power of $\phi$.
We define a new dimension for a subvariety $U$ of $V^n$ by $\dim_{zar}(U)=\dim(V) \inv \dim(U)$, and show by induction on $n$ that it is integral.  The dimension theorem follows
from the above property of intersections.}
    Even if $\dim(V)>1$, this derived structure will be intrinsically one-dimensional:
  it has no infinite, co-infinite definable sets.     From this it follows
 by a general theory that either modularity holds, or else a field can be found essentially
 as a quotient of a difference subvariety of $V^n$; in the latter case one shows that $\bV$
 is not field-free.  We refer to \cite{CHP} for details.  
 
 Note that taking all invariant subvarieties would not have worked:  the intersection of two
invariant $2$-dimensional subvarieties of a smooth $3$-dimensional one may be a union of two
$1$-dimensional parts, whose intersection in turn is a point.   If we allowed only invariant subvarieties, we would see two codimension $1$ sets intersecting in a codimension $3$ set,
contradicting a basic property of dimension.  The device of finding more regular behavior
by replacing $\phi$ by $\phi^n$ will be used repeatedly in the present paper; we will also
check that field-internality and modularity are insensitive to this.

  Let us relate modularity to some  degree-based notions used commonly in dynamics. 
 
\<{prop} Let $\bV = (V,\phi)$ be primitive. 
If $\deg(\phi)>1$ then   $\bV$ is fixed-field free.
If $\phi$ has separable degree $>1$ then $\bV$ is  modular.   \>{prop}
%
%\<{prop}  \<{itemize}\item[(1)]  If $\dim(V)=1$, then $(V,\phi)$ is fixed-field free if and only if 
%$\deg(\phi)>1$.  It is modular (or field-free) iff $\phi$ has separable degree $>1$.
%%Hence  in characteristic $0$,   $(V,\phi)$ is modular iff $\deg(\phi)>1$; for $V$ defined over a finite field, $(V,\phi)$ is modular iff $(V,\phi^n \circ Frob)$ has degree $>1$
%%for every $n$ and every Frobenius power Frob.  
%\item[(2)]
%   More generally, if $\bV$ is primitive and  $\deg(\phi)>1$ then $(V,\phi)$ is fixed-field-free.
%It is modular (or field-free) if $\phi$ has separable degree $>1$.
%   \>{itemize} \>{prop}
 If $V$ is a curve, then $\bV$ is primitive, and the converse of both statements is
 true.    But this is 
 a purely one-dimensional phenomenon.  In fact, if   $(V,\phi^n)$ is primitive for all $n$, and remains primitive after base change,   then $\bV$ is  {\em necessarily} modular when $\dim(V)>1$.    Thus for such strongly primitive $\bV$, modularity is equivalent to having either
 dimension or separable degree $>1$.  
  
If $V$ is not primitive, neither implication is true;  the condition $\deg(\phi)>1$ becomes rather weak, implying only 
the existence of one fixed-field-free subquotient in a decomposition.   
The very strong condition of ``polarization''  (cf.  \cite{CL})    implies that every subquotient
has positive degree and hence fixed-field-freeness, thus modularity
in characteristic $0$.

%
%What is new in the current presentation is   the fact
%that the theory - notably the definition and basic properties of internality - can be carried out at the level of rational maps rather than algebraic or definable closure.   The class of difference
%field extensions of $K$ is not an amalgamation class - $AD_K$ does not always admit products
%or fibered products - since such products may be reducible.  The usual approach is to
%form a theory without quantifier-elimination, and then analyze the definable sets.  This
%is what we do in \S 3, showing a posteriori that modularity and related notions do not
%depend on the choices made in the creation of a universal domain.  In an appendix \S 6
%we take an alternative approach, defining these notions directly for  appropriate expansions of stable theories, without assuming amalgamation.  A certain weakening of amalgamation does
%hold and is essential in both cases.   

A natural question is to what extent the  results described here can be formulated at the level of varieties, rather than birationally.   If a primitive difference variety is field-internal in the generic sense used here,
then in fact it contains a difference subvariety that is field-internal.  We do not know if the same holds for modularity.
%:e

 \para{\noindent \bf Isotriviality} 
 Let $K \leq L$ be difference fields, and let $V \in AD_L$.  
   Roughly speaking, we say that $V$ is isotrivial if it derives from an algebraic dynamics over $K$.  This splits into a number
 of technically distinct notions.  

\<{itemize}
\item[1.]
  $\bV $ {\em descends to $K$}   if it is $AD_L$- isomorphic to $\bW _L$
for some $\bW \in AD_K$.  

This is the finest notion we will use; as above  our isomorphisms are birational. 

\item[2.]  $\bV  $  is {\em constructibly isotrivial} if it is isomorphic to
$\bW_L$    using constructible maps, i.e. compositions of rational maps with inverses
of   purely inseparable rational maps.  

\item[3.]  %Let $\hK$ be a difference field. 
%  By a {\em finite cover} of $\bV$ in $AD_\hK$ 
%we mean an element of $AD_\hK$ admitting a generically finite morphism to $\bV$. 
  $\bV,\bV'$ are {\em isogenous} if they admit a common finite cover.   
$\bV  $  is {\em isogeny isotrivial} if there exists $\bW \in AD_K$
such that $\bV, \bW_{L}$ are isogenous in $AD_L$.  
 %(When $\bV/bW$  is the model-theoretically natural notion of non-almost orthogonality 
\>{itemize}
 
  In dimension one the notions (1),  (2) coincide, since a constructible bijection
 can always be made birational  by a Frobenius twist.   Beyond this remark
 we will consider only (2) and (3) here.   (3) is of course  too coarse to be of interest in pure algebraic geometry, but is the most basic equivalence relation to consider for varieties with a dynamics.

\<{prop} Let  $K  \leq L$ be algebraically closed fields, $\bV \in AD_L$.
%  Assume $V$ is absolutely irreducible, and
Assume no positive-dimensional quotient of $\bV$ has constant dynamics. 
%  Assume  $U,V$ are absolutely irreducible, 
%and $\bU=(U,Id)$ has periodic dynamics,
%while no 
%positive-dimensional quotient of $\bV$ does.  
Then   isogeny isotriviality is equivalent to constructible isotriviality for $\bV$.
\>{prop}
 
We sketch the proof. Assume condition (3) holds.  Then   there exists a common finite cover $\bY$ of $\bV, \bW_{L}$.  It follows that 
 $\bW_{L}$ has no positive-dimensional quotients with constant dynamics
 (cf. discussion following Lemma \ref{mod1}), hence the same is true for  $\bW$.
%$\bW=(W,\phi_{\bW})$. 
 By $(2^+)$,  there
exists a finite  cover $\bW'$ of $\bW$, 
with $\bW'_{L}$ a   finite cover of $\bY$.  Let $f: W' \to V$ be the $L$-definable map
coming from the composition $\bW'_{L} \to \bY  \to \bV$.  
Define $E$ on $W'$ by $wEw'$ iff   $f(w)=f(w')$.  This is a constructible equivalence
relation on $W'$ defined over $L$.   
Using the orthogonality \ref{mod1} (2)  % for $(W,\phi_W^n)$ 
we see that $E$ must be $K$-definable.  At the constructible level
this gives immediately a dynamics $\bW' / E \in AD_K$, isomorphic over $L$ to $\bV$,
hence constructible isotriviality (2).

 Given $\bU \in AD_K$ and $\bV \in AD_\bU$, we say $\bV$ is isotrivial (in any of these senses)
 if it is isotrivial in $AD_L$, with $L$ the function field of $\bU$.
The proposition applies   in particular when $\bV/\bU$ is modular, and  $\bU$ has constant dynamics.
 
\para{\noindent \bf Modularity and descent}  

Let  $\bU, \bV'' \in AD_K$, $\bV,\bV' \in AD_\bU$.    Say $\bV$ is {\em dominated by} $\bV' \in AD_\bU$ if there
exists a (dominant) morphism $\bV' \to \bV$ in $AD_\bU$; and by $\bV'' \in AD_K$ if
it is dominated by $\bV''_\bU$.  

\<{prop}  \vlabel{1.4}  Let $K=K^{alg}$, $\bU \in AD_K$   and  let  $\bW \in AD_{\bU}$ be  field-free.  
Assume $\bW$  is dominated by an object of $AD_K$.  Then $\bW$
is isogeny isotrivial.   \>{prop}
%
%  It  suffices to assume
%that $(W,\phi)$ is fixed-field free, i.e.   orthogonality to the Frobenius-twisted fixed fields is not needed. 
 
\prf   Say $\bW$ is dominated by $f: (\bV \times \bU) \to \bW$. 
Consider the dominant rational map $f: V \times U \to W$, ignoring the dynamics.
The following two statements are basic in Shelah's stability theory, and easy to
establish directly in the case of algebraic varieties.

(i) There exists a dominant $g: V \to V'$ such that 
$g(v)=g(v')$ iff  $f(v,u)=f(v',u)$ for generic $u \in U$ (i.e. whenever defined). The map
$f$ factors through $g$
and some $f': V' \times U \to W$.      For $v,v' \in V'$ we have $f'(v,u)=f'(v',u)$ for generic $u \in U$ iff $v=v'$.    

(ii)  Define $f_m: V' \times U^n \to W^n$ by $f_m(v,u_1,\ldots,u_n) = (f(v,u_1),\ldots,f(v,u_m))$.  Then  for some $m$, $(f_m,Id)$ embeds $\bV' \times \bU^m$  into $\bW^m$ (generically).   (Follows from the last statement of (i)). 

Now recall $\phi_V,\phi_W$.  If $g(v)=g(v')$ then $g(\phi_V(v)) = g(\phi_V(v'))$:  indeed 
if $f(v,u)=f(v',u)$ for generic $u$, then
$f(\phi_V(v),\phi_U(u)) = \phi_W( f(v,u) ) = \phi_W( f(v',u)) = f(\phi_V(v'),\phi_W(u))$;
and $\phi_W$ was assumed dominant.   So we can define a rational $\phi': V' \to V'$ with $\phi'(g(v))=g(\phi_V(v))$.   Let $\bV' = (V',\phi')$.    By (ii) we have $f(V) \leq \bW \times _\bU \cdots \times_\bU \bW$.  By  Proposition \ref{mod1} (1) over $\bU$, $\bV'$ is field-free, hence modular.  

  Now view $\bV'$
as parameterizing a family of functions $\bU \to \bW$; their graphs are irreducible subvarieties
of $\bU \times \bW$ of dimension equal to $\bU$; so by modularity, $\bV'$ may again be replaced
by a quotient $\bV''$ of dimension $\dim(U \times W) - \dim(U)= \dim(W)$.  Hence
the graph of $f$ is a finite cover of $\bV''_\bU$ and also of $\bW$, showing they are isogenous, and isotriviality (3) holds.     \qed

Each of the steps in this proof is valid in great generality, in categories of definable sets
of theories such as ACFA.    Given the basic properties of such definable sets, 
the proof generalizes to  dynamics given by correspondences.    This will be done in 
Part II.

%An   model-theoretic property that will not be discussed here is {\em stability}.  Let us just
%mention that for algebraic dynamics, it is equivalent to the non-existence of moving families
%of finite covers in the $AD$ category. In characteristic $0$, it is shown in \cite{CH}
%to be equivalent to modularity.  At all events when $\bW/\bU$ is stable, then 

  {\em Assume now 
that the dynamics on $\bU$ is trivial.  Then in Proposition \ref{1.4} one can conclude
that $\bW$   descends constructibly to $K$. }   We have a morphism $f: \bU \times \bV \to \bW$ in  $AD_\bU$, and as above we may take $\bV$ to be  fixed-field-free. Given $u \in U$ we have a constructible equivalence
 relation $E_u$ on $\bV^2$, namely $xE_u y$ if $f(u,x)=f(u,y)$.  If $u \in \bU(\Uu)$
 then $E_u$ is compatible with $\phi_V$, as we saw in the proof above.  By   the orthogonality
 principle of Proposition \ref{mod1}, $E_u$ cannot really depend on $u \in \bU(\Uu)$; so $E_u=E_{u'}$ for generic $u,u' \in \bU(\Uu)$, and hence $E_u=E$ generically for some $E$.  Now   $\bW$ descends constructibly to $\bV'=\bV/E$.
   
   The same argument shows that $\bW$ descends together with any additional structure
   it may carry; notably if $\bW=( W,\phi)$ and $\phi': W \to W$ commutes with $\phi$,
   then $(W,\phi,\phi')$ descends.  For the image of the graph of $\phi'$ under
   the isomorphism $\bW \to \bV'_{\bU}$ is a  difference subvariety
of $\bV^2$, definable over the generic point of $\bU$;   as noted above for $E_u$, any such variety is $K$-definable.    Thus it suffices to  assume that $ (W,\phi^l)$
is dominated by an object of $AD_K$, for some $l$ (for then $(W,\phi^l)$ descends, and
$\phi$ commutes with $\phi^l$.)

Part II contains more general statements.   
 We show by example  (3.5 in Part II) that Proposition \ref{1.4} (or \ref{1.9} below) fails without a modularity or orthogonality assumption, even when 
the dynamics on $\bU$ is trivial.  

 \para{\bf \noindent Dynamical Northcott  for function fields.}
 The assumption of domination by a difference variety over $K$ can be rephrased 
 in the language of limited subsets of the function field
 $L=K(U)$.  If $g: V \times U \to W$ is the dominating map, and $a$ a generic point
 of $U$, then the image of $V(K)$ under  $x \mapsto g(a,x)$ is a typical {\em limited subset of
 $K(U)$.}   When $V(K)$ has an Zariski dense $\phi$-orbit, $W(L)$ will have a Zariski dense
 orbit contained in a limited set.  Conversely, if $W(L)$ will have a Zariski dense
 orbit contained in a limited set, then $W$ is dominated by a difference variety over $K$.  
 We thus have:

{\bf \noindent Corollary}  \vlabel{1.9}  Let $(W,\phi)$ be a fixed-field-free   algebraic dynamics over $L$, $L$ a finitely generated extension field of  $K=K^{alg}$.   Assume $(W,\phi)$ does not constructibly descend to $K$.
  Then no limited subset of $W(L)$ can  contain a   Zariski dense $\phi$-orbit.

 The case $V=(\Pp^1)_L$,   is the result  of \cite{baker}.   In fact the 
 hypothesis $K=K^{alg}$ is not needed when $L/K$ is regular, answering a question posed there for $\Pp^1$.

\<{thm} \vlabel{g2}  Let $(V,\phi)$ be a primitive algebraic dynamics over $L$, $L$ a finitely generated 
regular extension field
 of a field $K$.  Assume %$V$ is absolutely irreducible and reduced, and 
 $(V,\phi)$ does not constructibly descend  to $K$.  Then no limited subset $Y$ of $V(L)$
 contains a  Zariski dense $\phi$ orbit; in fact for some $n=n(Y)$  and a finite number
 of   proper subvarieties $U_1,\ldots,U_j$ of $V$, defined over $L$,
 there is no $a \in V(L) \m \union U_i(L)$
 with $a,\phi(a),\ldots,\phi^n(a ) \in Y $. 
  \>{thm}
  
  The finite bounds  $n,j$ follow  by compactness.  In some cases,  there exists canonical height, i.e. a function $h: V(L) \to \Rr^{\geq 0}$
 such that the inverse image of a compact set is limited, and $h(\phi(a)) = \kappa h(a)$
 for some $\kappa >1$.  In this case, taking $Y = h 
 \inv ([0,1]) $, there will be no $a \in V(L) \m \union U_i(L)$ with
 $h(a)< \kappa^{-n(Y)}$.   
 
We first prove this theorem under  the additional assumption that 
 $(V,\phi)$ is fixed-field-free.  
The methods presented in  the present part suffice to prove this   for isogeny-isotriviality.  Using Theorem 3.3 of Part II, 
we improve to constructible isotriviality; see \S  \ref{limited}  \footnote{Theorem 3.3 of Part II assumes $\deg(\phi)>1$, but in fact uses this only via the conclusion of  Proposition 1.4, that $\bV$ is fixed-field-free.  It is in this form that we use it here.}.  
Finally in Part III,   the primitive,   fixed field internal dynamics are   explicitly described;
  they are associated  with translations in one-dimensional algebraic groups; 
  using the classification of these, it turns out that the fixed-field-free assumption 
  can be removed.   
  
%We discuss these   statements in   \S  \ref{limited}. 

\para{\noindent \bf Difference varieties}
 
If $(V,\phi)$ is an algebraic dynamics over $K$, the function field $L$ of $V$ is a difference field 
that is finitely generated over $K$ {\em as a field}.  
This can be made geometric in one of two ways:
\begin{itemize}
\item[1.]  
  One can consider pairs $(V,\phi)$
with $V$ a pro-algebraic variety.  
\item[2.]  one can stay with varieties but replace the rational map $\phi$
by a   {\em correspondence }  on $V$, i.e. a subvariety $S$  of $V \times V$
with generically finite projections.  For any difference field 
$(L,\si)$ extending $K$ it is possible, using the Ritt-Raudenbusch
basis theorem,  to find $(V,S)$ such that $L=K(a)_\si=_{def}K(a,\si(a),\si^2(a),\ldots)$ 
 with  $(a_0,a_1)$ a generic point of $S$,  and such that moreover any other difference field  
of this description is isomorphic over $K$ to $L$.   
The extension to such ``nondeterministic dynamics'' should allow for greater flexibility.
\>{itemize}

All our structural results remain valid in the more
general context, when $L$ is finitely generated over $K$ as a difference field.
Theorem \ref{g2} applies to difference varieties presented via a correspondence
$(V,S)$, but in the following form: for some $n$, and some finite  union $Z$ of proper subvarieties, there is no $n$-cycle $a,\phi(a),\ldots,\phi^n(a) $ each of which is
in $Y \meet V(L)  \setminus Z$.  In other words, for some $Z'$ and any $a \in V(L) \setminus Z'$, no  choice of $\phi$-preimages can be iterated $n$ times.

%In \S 2 we will recall the basic theory of difference fields relevant to this paper.  \S 3 includes
%the necessary work to move to rational maps, concluding in particular with 
%  Theorem \ref{mod} as stated above, and variants.  \S 4 draws corollaries on descent.
%The expository \S 5 develops the language of limited sets, leading to a restatements of the
%descent results, and in particular generalizations of Baker's theorem.  

\para{\bf \noindent Definable Galois theory}
Some of our results have been relegated to a future paper.  This includes
more precise descriptions of non-disintegrated modular sets, and of 
 of field-internal dynamics.   We briefly discuss the latter here.   The issue is the elaboration in the present context of definable Galois theory, introduced into model theory, in the differential setting, by Poizat.
 
%
%Let $K$ be  a difference field, 
%$G \times X \to X$ be an algebraic group action defined over  the field of constants $C_K$ of $K$.
%%, and $\bX$ a 
%  %homogeneous space for $\bG$, 
%  % a  $C_K$-variety with a transitive $\bG$-action. 
%Let
%$\fg \in \bG(K)$, and let $T(\fg)$ denote the function $x \mapsto g \cdot x$.  Then
%the difference dynamics $(X,T(\fg))$ is field-internal.   This will be called a {\em translation variety}.

Let $K$ be a field, and let  $\bV \in AD_K$ be field-internal.  Assume $V$ is absolutely
irreducible.  Then   $\bV$ can be shown to be   a translation variety, or fibered over
a periodic dynamics by translation varieties.  The automorphism
group of $\bV(\Uu)$ over the fixed field of $\Uu$ can be seen to be a part of 
an algebraic group $G$, and $\phi$ will be an element of (another part of) $G$.

In the relative case, when $K$ is a difference field and $V$ is not necessarily defined
over the fixed field, more definable Galois theory is needed, including the theory of the opposite
group (a group over the constants, isomorphic after base change to the automorphism group) and the existence of a bi-torsor for the two groups.  We obtain essentially the same result, reducing field-internal dynamics to translation varieties, when the 
fixed field  of $K$ 
is pseudo-finite. The pseudo-finiteness condition is analogous to the requirement in differential
Picard-Vessiot theory that the field of constants be algebraically closed.  
Over more general fields, delicate questions of Galois cohomology arise.

\bigskip

	%
	%\section{From types to isomorphism types}

	%As defined in \cite{CH},\cite{CHP}, modularity and fixed-field internality were properties of
	%a complete type  $tp(c/K)$ in a saturated model   $\calu$ of ACFA.  We show here that they actually depend only on the difference field extension $K(c)/K$.  In the Appendix, we present
	%the same material differently, defining the notions directly in a way that does not use the embedding.  

	%
	%\def\proj{\Pp}
	%\noindent{\bf Example}    An algebraic dynamics over $k(t)$, dominated by an
	%algebraic dynamics over
	%$k$, but not isogenous to a difference variety over $k$.    \rm

	%Let $H$ be a vector extension of an Abelian variety; i.e. there
	%exists an exact sequence of algebraic groups $0 \to V \to H \to A \to 0$ with $A$ an Abelian 
	%variety, $V \cong G_a^n$ a vector group.  Assume $\dim(V) = 2$ and 
	%$Hom(H,G_a) = (0)$.  Assume $H$ is defined over $k$, and identify the projective
	%space $\proj V $ with  $\proj^1$.   In particular a transcendental element $t$ gives a one-dimensional subspace $V_t$ of $V$.   Let $H_t = H / V_t$.  Fix $h \in H$ generating $H$;
	%it suffices that the image of $h$ in $A$ generate $A$; if $A$ is a simple Abelian variety, it
	%suffices therefore that the image is nonzero.  Let $Y=  (H,T(h))$ and $X_t = (H_t, T(h_t))$
	%where $h_t$ is the image of $h$ in $H_t$, and $T(g)$ denotes translation by $g$.  Then $Y$ dominates $X_t$.  But $X_t$
	%is clearly not isotrivial.   % \>{example}
  
%In positive characteristic, we could similarly factor out  a finite subgroup of $V_t$, giving 

\bigskip\noindent

\section{Internality without quantifiers} 

   In the context of a universal domain $\Uu$, one says
that a type $P$ is {\em internal} to a definable set (or union of type-definable sets)
$\pi$ if any realization of $P$ is definable over a fixed finite set, and elements of $\pi$.  
The theory began in the stable setting (\cite{h-s}), with three   components:

a)  Existence of a canonical maximal quotient $P_\pi$ of $P$,  internal to $\pi$.

b)  Domination:  any interaction between a realization $b$ of $P$
and elements of $\pi$ must involve   the image of $b$ in $P_\pi$.  

c)  Definable Galois theory:   the interaction of $P_\pi$ and $\pi$ is  controlled
by $\infty$-definable automorphism groups, or in a more precise version, groupoids.  

Parts (a) and (b) had their origins in Shelah's semi-regular types (based on acl rather than dcl).
The theory was generalized to arbitrary first order theories in  \cite{h-b} (appendix B), and
part (c) was transformed to  a quantifier-free setting in    \cite{kam2}; but a general treatment
of (a,b) is without quantifiers is missing. 

In the simple context, all parts of the theory were generalized, at least ``up to algebraic (or bounded) closure"; see \cite{HKP}, \cite{SW}, \cite{BW}.   Roughly speaking, (a,b)
are carried out for a certain notion of internality,  (c) for a finer one, and
the internal parts $P_\pi$ in the two senses are  shown to be related by multi-valued maps.  
This corresponds precisely to what we do in Parts II (q.f.-internality) and III (definable
Galois theory) , except that in our context
we are able to identify the finer version as quantifier-free internal; for our purposes it is vital to identify the maps as rational maps.   

Simple theories are often obtained from expansions of stable theories by an amalgamation
process; the additional structure results from the failure of amalgamation over non-algebraically closed sets.  We place ourselves in such an enriched stable context here.  In place of amalgamating, defining internality using quantified formulas, and then   show ing that  it depends only on the quantifier-free part, we shortcut and define internality before amalgamation.   Examples include any number of derivations and automorphisms, commuting or not.

Simplicity corresponds essentially to uniqueness of amalgamation of the expansions,
relative to uniqueness for the stable layer.  This is valid in the examples mentioned above.  We do not assume this, except for
Remark \ref{!}.   Our framework thus  includes non-simple structures, but we have not investigated what it means for them.

We consider below a class $\mC$ of $\mL$-substructures,   with a stable
reduct to a language $\mL_0$.   We assume $\mC$ is closed under substructures.  The stable reduct 
is used  to obtain a notion of free amalgamation on $\mC$, enjoying the usual
properties, and in particular the existence of canonical bases:  for any structure $C$
in the class and any $A,B \leq C$, there exists a unique minimal substructure $B'$ of $B$
such that $A,B$ are independent over $B'$.   It would be possible to replace the ``stable sublanguage" hypothesis by assuming abstractly that $\mC$ is given with
an amalgamation notion, having the usual properties, and with canonical bases.

\para\vlabel{int11}{\bf Internality and orthogonality}.

Let $\mT_0$ be a stable theory in a language $\mL_0$. 
We assume $\mT_0$ eliminates quantifiers and imaginaries, and that substructures of models of $\mT_0$  are 
 definably closed.  
 
  We recall some facts and notation:  let $K,L,L'$ be 
 substructures of some model of $\mT_0$.    
 The substructure generated by $K \union L$ is denoted $KL$.  
  We write $L^{alg}$ for the algebraic closure of $L$.
  We say 
 $K_1,K_2$ are strongly free over $L$ if they are independent over $L$,
 and whenever $K_1',K_2'$ are independent over $L$ and $\a_i: K_i \to K_i'$
 is an $L$-isomorphism, then $\a_1 \union \a_2$ is an $L$-isomorphism. We say
 $L/K$ is {\em stationary}
if $L$ is strongly free from $K^{alg}$ over $K$; equivalenly,
$\dcl(L) \meet \acl(K) = \dcl(K)$. 

$Cb(K/L)$ denotes the smallest substructure $L'$ of $L$ such that $K, L$ are strongly
free over $L'$; equivalently they are free over $L'$, and  $ \dcl(L',K) \meet \dcl(L) \subseteq \dcl(L')$.  
If $e \in Cb(K/ (L')^{alg})$ then $e \in KL$; let $f$ code the  the finite  set of realizations of
$tp(e/L')$;  then $f \in Cb(K/L)$; and such elements $f$ generate $Cb(K/L)$.  Equivalent forms of these notions appear in   Shelah and in 
Bouscaren \cite{bouscaren}.  

  Let $\mL$ be a language containing $\mL_0$.  Let $\mC$ be a class of $\mL$-structures, closed under isomorphism and substructures.

The reduct of an $\mL$- structure $L$ to $\mL_0$
is denoted $L_0$.     Let $\mC_0=\{L_0: L \in \mC \}$.  We assume for the rest of this section:

(I0)  If $A \leq B,C \in \mC$, and $B_0,C_0$ are strongly free over $A_0$, then their free
amalgam $D_0$ can be expanded to some $D \in \mC$, with $B,C \leq D$.    Conversely, for any 
$B,C \leq D \in \mC$ we have $(BC)_0=B_0C_0$.

Let us immediately explain how (I0) will be used.  Assume $B^1,B^2,B^3, \ldots \in \mC$
and $B^1_0,B^2_0,\ldots$ are free over $A_0$.  Let $A_0^i = A_0^{alg} \meet B^i_0$,
and let $A^i$ be the substructure of $B^i$ with reduct $A_0^i$.  
Assume there exists isomorphisms $f^i: B^1_0 \to B^j_0$ such that $f^i | A^1$ is 
an isomorphism $A^1 \to A^i$.   So  all $A^i$ are copies of a single structure $A'$.  
Now $B^1_0,B^2_0,\ldots$ are strongly free over $A'_0$.
By (I0), their free amalgam $D_0$ can be expanded to some $D \in \mC$,
with each $B^i \leq D$.    In particular, the $B^i$ can be amalgamated over $A$
within $\mC$.  If only $B^1$ is given, it is always possible to 
find $B^2_0,\ldots$ free over $A_0$, and expand them to $B^2,\ldots$
in such a way that even $B^1,B^i$ are isomorphic; so $B^1,B^2,\ldots$
are jointly freely embeddable into an element of $\mC$.

We will write $K,L,K',J$ etc. for elements of $\mC$.   We let $\qfcb(K/L)$
be the $\mL$-structure generated by $Cb(K_0/L_0)$.  $tp_0(a/L)$ denotes the $\mL_0$-type.

Let  $\pi$ be a set  of $\mC$-isomorphism types.   We write $L \in \pi$ to mean that the isomorphism
type of $L$ is in $\pi$.  We assume: 
 
(I1)  Let $L \in \mC$ be  generated by some collection ${\mathcal A}$ of substructures of $L$.
  Then $L \in \pi$ iff ${\mathcal A} \subseteq \pi$.   
  
If $L \in \mC$, let $L_\pi$ be the join of all 1-generated substructures of $L$ in $\pi$;    this is the largest substructure of $L$ in $\pi$.

Consider the properties:

(I2)  Let 
%M
 $K,J \leq L   \in \mC$, with $J \in \pi$ and  $J_0 \leq K_0^{alg}$.      Then there exists    $N \in \mC$, $J \leq N$, with $N_0 \leq K_0^{alg}$, and such that $N_\pi$ is $Aut(K_0^{alg}/K_0)$-invariant. 

(I3)   Let $L \in \pi$ and let $L' \in \mC$.  Then $\qfcb(L/L') \in \pi$. 

\claim{}    If   $\pi$ has (I1),(I2), it also enjoys (I3). 

Indeed let $L_0,L^1_0,\ldots$ be an independent sequence over $L'_0$ in some model   of 
$\mT_0$, such that all $L^j_0$ satisfy the same type over $(L'_0)^{alg}$.
 Expand each $L^i_0$ to $L^i$ so that $L^i \cong L$.  Then each $L^i \in \pi$,
and by  the remark following (I0)  they jointly embed  into some $L^*_0$, $L^* \in \mC$.    Let
 $E=Cb(L_0 / (L'_0)^{alg}) \subseteq  L^1_0L^2_0 \cdots$.   By  (I1) (the `if' part) we have $L^1L^2 \cdots \in \pi$.   Since $\pi$ is closed under substructures (the `only if' part of (I1)), we have  $E \in \pi$.  
For $e \in E$ we have $e \in  (L'_0)^{alg}$; let $e' \in L'_0$ code the $Aut( (L'_0)^{alg} / L'_0)$-orbit of $e$; 
by (I2) we have $e' \in L^*_\pi$.         By the discussion of canonical bases
 in stable theories above,  the set of such $e'$ generates $Cb(L_0/L'_0)=\qfcb(L/L')$.  

We will say that $A,B$ are (strongly) free over $C$ if $A_0,B_0$ are so.

\<{defn} \vlabel{def-int}Let  $J \in \mC$.     
  \<{itemize} \item  
We say that $J/K$ is   {\it qf-internal to $\pi$} if
there   exists $\hL $ with $\hL_0 / {J_0}$ stationary,  and   
  $K \leq L \leq  \hL$,  such that $J, L$
are free over $K$,   
and $\hL =L \hL_\pi =L J$.   
\item 
Let $K \leq K_1 \leq K_2$.  We say that $K_2$ is $\pi$-dominated by $K_1/K$ if
for any $\hL \in \mC$   with $K_2 \leq \hL$ and $K \leq L \leq    \hL$, 
  if $K_1$ is free from $L$ over $K$ then   $K_2$ is  strongly free  from 
  $K_1 L \hL_\pi$ over   $K_1L$.  \>{itemize} \>{defn}

 \<{thm}  \vlabel{6.3} Let $K \leq J \in \mC$.  Then there exists a unique sub-extension  $K_1$ with 
 with $K_1/K$ qf-internal to $\pi$, and $J$ $\pi$-dominated by
 $K_1/K$.     
\>{thm}  

\prf   Let $K_1$ be a sub-extension of $J/K$.   Call $K_1/K$ {\em  weakly qf-internal to $\pi$}  if
there exist  $K \leq L,K_1 \leq  \hL_1$ with  $\hL_1/K_1$  
%M
 stationary,  $K_1, L$
  free over $K$,   and $K_1 \subseteq L (\hL_1)_\pi $.   

\claim{1}  If $K_1/K$ is weakly qf-internal, and $J$ is $\pi$-dominated
by $K_2/K$, then $K_1 \subseteq K_2$.   

To see this,  let
$K \leq L,K_1 \leq  \hL_1$ be such that  $\hL_1/K_1$ is stationary,  $K_1, L$
are free over $K$, and and $\hL_1 =L (\hL_1)_\pi$ contains $K_1$.  
We may jointly embed $\hL_1,J$
into some $\hL \in \mC$ with $K_1 \leq \hL$ and $\hL_1,J$ free over $K_1$; hence $L,J$ are free 
over $K$; we may take  $\hL = J \hL_1$, so $\hL  = J L (\hL_1)_\pi   $ and
in particular $K_1 \subseteq L \hL_\pi$.

By definition of domination, $J$ is strongly free over $K_2L$ from $K_2L \hL_\pi$.
As $K_1 \subseteq L \hL_\pi$ we have $K_1 \subseteq K_2L$.  But $K_2 \leq J$ and (by stationarity) $J,L$ are strongly free over $K_1$ so $K_1 \subseteq K_2$.  
 
Uniqueness follows immediately from Claim 1,  since if $K_1,K_2$ are two candidates,
then $K_i/K$ is weakly qf-internal to $\pi$ and $J$ is $\pi$-dominated by $K_{2-i}/K$,
so $K_i \subseteq K_{2-i}$.

Existence:  Let $K_1$ be the join of all  sub-extensions of $J/K$ that are 
qf-internal to $\pi$.  By amalgamating the relevant witnesses $L$ over $K$,we  see
that $K_1/K$ is qf-internal to $\pi$.  We have to show that $J$ is $\pi$-dominated 
by $K_1/K$.   In other words, consider $J,L \leq \hL$, with $L$   free from $K_1$ over $K$.
We have to show that  $M=\hL_\pi$  is strongly free from  $J$ over $LK_1$.

Suppose otherwise.  Then there exist   $J,L \leq \hL$, with $L$   free from $K_1$ over $K$,  and $M=\hL_\pi$ not strongly free from  $J$  over $LK_1$.   
 
 \claim{2}  We may choose $L,\hL$ so that $L J /J$ is stationary.
 
 \prf      For $i=1,2,\ldots$ let $J^i,M^i \in \mC$, contained in some $L^* \in \mC$ with $L \leq L^*$ and  $J^i M^i \cong_{L^{alg}} J M$ by an isomorphism
 taking $J , M$ to $J^i,M^i$, and with $J^iM^i$ an independent  sequence of extensions of $L$.  Take $J^0=J, M^0=M$.  Then the  $J^i$  are  independent over $L$, and     each one is independent from $L$ over $K$, so $L,J,J^1,\ldots, $ are  independent over $K$.    Let $L' = J^1  J^2   \ldots$. 
  Since the $J^i$ have the same type over $K^{alg}$,
 it follows that $J^i /  J$ is stationary; hence $L' / J$ is stationary.  
 Let $E=\qfcb(J M/L)$. 
 So $J M, L$ are strongly free over $E$.  Thus $J,ME$ are not strongly free over
 $K_1$, hence (since $J,L'$ are strongly free over $K_1$) 
 $J,ME$ are not strongly free over $L'K_1$. We have 
  $E \subseteq J^1M^1J^2 M^2 \cdots = L'M^1M^2 \cdots $,
so $ME \subseteq L' M^*_\pi$.  Thus  $J, L' M^*_\pi$  are not strongly free over $L'K_1$.   \eprf

Assume now that $L,\hL$ are as in Claim 2. Let $M' = Cb( M / L  J)$  Then by Claim 2,
$M'$ is not contained in $LK_1$.  On the other hand $M' \subseteq L J$ and
 by (I3), $M'/K \in \pi$.   
 
 Finally let $K_1' = Cb(M'L / J)$.  Then $K_1'/ K $ is qf-internal to $\pi$.  For $L/K$ is 
 stationary, and if $L_i,M'_i$ are indiscernible independent copies of $M'L/J$,
 $L_\infty = L_1L_2 \cdots, M'_\infty = M'_1 M'_2 \cdots$ then
 $L_\infty M'_\infty = L_\infty K_1'$.  By maximality of $K_1$ we have $K_1'=K_1$.  
 Since $M' \subseteq L J$ we have $M' \subseteq LK_1$.  This contradiction proves the proposition.  

\eprf

Here is a corollary to the proof.  (2) is a version of the ``Shelah reflection principle"
referred to in the introduction.

\<{cor} \vlabel{c12} Assume (I1,I2).  

\noindent (1)  weakly qf-internal is the same as qf-internal.

\noindent (2)  Let $J/K$ be finitely generated.  In the definition of qf-internal (\ref{def-int}), one can take $L$ to be 
the amalgam of finitely many freely joined indiscernible copies of $J/K$.
 
\>{cor}

It will be worthwhile to review the statement of (2) in terms of definable   maps.
 Since 
$J_0 / K_0$ may not be finitely generated, we need to work in the pro-definable
category (cf. \cite{kam}); this is really a side issue here, and the reader may prefer at first to assume
finite generation or ignore it.
By an $\mL_0$-morphism
(or an $\mL$-morphism) we mean a term of $\mL_0$ (respectively $\mL$).
    Assume the language $\mL_0$ has constants for the elements of $K_0$. 
Let $\mC_0^*$ be the category whose objects are complete quantifier-free *-types of
$\mL_0$.   We can think of an object of $\mC_0^*$ as a pair $(L,a)$
with   $L \in \mC$ and $a$ a sequence of elements of $L_0$, generating $L_0$, up to isomorphisms
$(L,a) \to (L',a')$ with $a \mapsto a'$.    A {\em morphism} from $(L,(a_i)_{i \in I}) \to
(L,(b_\ups)_{\ups \in \Upsilon})$ consists of a choice of functions $f_\ups$
given by terms in $\mL_0$, such that $f_\ups=f_\ups( (x_i)_{i \in I(\ups)})$ is defined on $(a_i)_{i \in I(\ups)}$
and $f_\ups((a_i)_{i \in I(\ups)}) = b_\ups$.   Given $p_0 \in \mC_0^*$ we have the usual notion of
the type of an $n$-tuple from an indiscernible, independent sequence in $p_0$; it is denoted
$p_0^{(n)}$, and is canonically determined even if $p$ is not stationary.  Here is the restatement of Corollary \ref{c12} (2).   Below, as usual, $K(c)$ denotes 
 the substructure of $J$ generated by $c$ over $K$.

\<{cor}  \vlabel{h}  $J/K$ is qf-internal to $\pi$ iff there exists $p=(J,a) \in \mC_0^*$,   
and an invertible morphism $h$ of $\mC_0^*$ with domain $p_0^{(n)}$, 
%$q= (\widehat{L},b,c) \in \mC_0^*$ with $(\widehat{L}_0,b) = p_0^{n-1}$ and $K(c) \in \pi$,
% p_0^{(n)} \to q_0$, for some
%$q= (\widehat{L},b,c)$ with $(\widehat{L}_0,b) = p_0^{n-1}$ and $K(c) \in \pi$; and 
such that $h(a_1,\ldots,a_n)=(a_1,\ldots,a_{n-1},c)$ with $K(c) \in \pi$.  
%: \<{itemize}
%\item $p'_0 = p_0^{(n)}$
%\item $q= (\widehat{L},b,c)$ with $b \models p^{n-1}$ and $K(c) \in \pi$.  
%\>{itemize}
\>{cor}

 This gives the assumption used in \cite{kam2} to obtain the quantifier-free liaison group.
 
 % The terms ``pro-definable" (Grothendieck) or *-definable (Shelah) mean:   defined by a set of formulas, with a possibly infinite number of variables.  A morphism  between pro-definable sets is given by an appropriate
%system of morphisms between definable factors; see \cite{kam}.    For instance
%a morphism from $\Pi_{i=1}^\infty D_i$ to $D'$ is just a term function from
%$D_1 $

%

%%Given a partial quantifier-free *- type $P$, let $L_P$ be the substructure of $L$ generated
%%by all realizations of $P$ in $L$.  
%Let $L(p_0;n)$ be the $\mL_0$-structure generated by  realization of
%$p_0^{(n)}$.   Let $\mC_0(p_0) = \{L(p_0;n): n=1,2,\ldots \}$. 
% Let $\mH(p_0;n)$
%be the set of morphisms  $h$ of $\mC_0^*$ with domain $p_0^{(n)}$, 
%with  $h(a_1,\ldots,a_n)=(a_1,\ldots,a_{n-1},c)$ with $K(c) \in \pi$; let $h_n$ be the last coordiante of $h$, i.e. 
%$h_n(a_1,\ldots,a_n)=c$.     
%$ be the class of all $L_0 \in \mC_0$ generated by a realization of
%$p_0^{(n)}$.  

%\claim{}  There exists a $\mT_0$- pro-definable group $\tG$ and a pro-definable action of $\tG$ on $p_0$, with the following property:   for any $n$ and any
%$\mL_0$-structure $L_0  \in L(p_0;n)$, $\tG(L_0) = G(L_0)$.  

% 
%Let $G(L_0)$ be the set of $\mL_0$-automorphisms $\alpha$ of $L_0$ such that
% if  $a=(a_1,\ldots,a_n) \models p_0^{(n)}$ is a generator of $L_0$,
%and $h$ is an $\mC_0^*$-morphism with $h(a)=(a_1,\ldots,a_{n-1},c)$,  
%then $h(\a(a_1),\ldots,\alpha(a_n))  = $\a(a_1),\ldots,\a(a_{n-1}),c$.  

%
%for any $h,a$ as in Corollary \ref{h}, $h_n(\a(a_1),\ldots,\alpha(a_n)) = h_n(a_1,\ldots,a_n)$.   

 \smallskip
 
 \para{\noindent \bf Internality in a universal domain}  
  Let $\calu$ be a universal domain for $\mC$, in the sense that if $A \leq \calu$ is small,
 $A \leq B \in \mC$, and $B_0/A_0$ is stationary then $B$ embeds into $\Uu$ over $A$. 
     We say $tp(a/K)$ is   qf-internal to $\pi$  if $K(a)/K$ is.  Let $D=D_\pi  (\Uu)= \{a: K(a) \in \pi \}$.  
We also speak of qf-internality to $D$.  

Say (I0!) holds if the expansion in (I0) is unique.  
 For instance this is the case if $\mC$ is the class of fields   with some distinguished derivations and automorphisms, see below.

\<{remark} \vlabel{!}  Assume (I0!, I1).  Then 
 qf-internality to $\pi$  implies that $tp(a/K)$ is $D$- internal in $\calu$, i.e.  for some fixed finite set $F$   of parameters, every realization of $tp(a/K)$ is definable over $F$ together with finitely many
elements of $D$.  \>{remark}
  
The theory of liaison groups for qf-internal types will be described in the sequel.
By a  $\mC$-group we will mean an
 $\mL$- quantifier-free definable subgroup $G$ of an $\mL_0$-definable
 group $H$.  We will show the existence of an $\mC$-group $G$ such that 
 for any $L \in \mC$, $G(L)$ is canonically isomorphic to the image of $Aut(L/ \pi(L))$
 in $Sym(p)$.  Moreover there exists a $\pi$-groupoid $H$ with torsor $Q$ such
 (fibered over an object set $Y$) such that $H_y$ and $G$ act  regularly on $Q_y$,
 by commuting actions.  The orbits of $G$ on $p$ do not always coincide 
 with the quantifier-free types over $\pi$, but if $T_0$ is $\omega$-stable they cut them
 in finitely may pieces.

\para{\bf }    We now check that (I1,I2) hold for difference-differential  fields.
 
 Let $\mT_0=ACF_F$ be the theory of algebraically closed fields over some base field $F$. Let $\mL$ be the language
 of fields with $l+l'$   operators  and let $\mC$ be a class of $\mL$- structures
 where the first $l$ operators are automorphisms and the rest are derivations.  We may ask that certain pairs commute.   A central case is that $\pi$ is the subfield fixed by some of the $l$ difference operators and annihilated by some of the $l'$ derivations.  If $l>1$ we will assume
 that we are in this case.  
 %should be able to do the commuting case in general.
 
   For $k \in \Nn, L \in \mC$, let   $L[k]$ denote   $L$ with $\si$ replaced by $\si^k$.

% Let $K$ be a difference field, $\mC$ the set of (isomorphism types of) difference field extensions of $K$. 
 
  Recall that $L/K$ is called primary if $L_0/K_0$ is primary, i.e. iff  $L_0/K_0$ is
stationary.   
 
\<{lem}  \vlabel{2-1}  Let $(L,\si)$ be any difference field, $K$ a difference subfield,
$k= K \meet Fix(\si)$, 
and $a \in Fix(\si) \meet K^{alg}$.  Then $a \in k^{alg}$.  \>{lem}

\prf Let $F \in K[X]$ be the monic minimal polynomial of $a$ over $K$.  Then $F^\si(a)=F^\si(\si(a))=
\si(F(a))=0$ and $F^\si$ is monic, so $F^\si =F$. Hence $F \in k[X]$.  \eprf

\<{lem} \vlabel{2-2'} Assume  (I1) holds for $\pi$.  Then  (I2) holds too assuming: 

\noindent (I2') \  Assume $L[k] \cong_{K[k]} L'[k]$.  Then $L \in \pi$ iff $L' \in \pi$.   \>{lem}

%M
\prf    Let $L,K,J$ be as in (I2); we may take $K=L(a),J=L(a')$.  
  Let $N$ be the normal closure of $L(a)$ over $L$, $\tau \in Aut(N/L)$ with $\tau(a)=a'$.  By  \cite{CHP} (1.11), for some $\ell$,
$\tau$ extends to an automorphism of the $\si^\ell$-difference field generated by $N$.
So $L(a)_\si[\ell] \cong _{L[\ell]} L(a')_\si[\ell]$.  In particular
  $K(a)_\si [\ell] \cong_{K[\ell]} K(a')_\si [\ell]$.  Since $K(a) \in \pi$, by (I2') 
  we have $K(a') \in \pi$.  
\eprf

 For several difference operators, one may consider an analogous condition,
considering subgroups of finite index in place of $\si^k$.  Does the lemma remain valid?
At all events, (I2) certainly holds if $\pi$ is the fixed field of some of the operators;
this follows directly from Lemma \ref{2-1}.

An element of a difference field is (twisted) {\em periodic} if $\si^m(a) = a^{p^r}$
for some $m \geq 1$ and $ r \in \Zz$, where $p=1$ in char. 0 and $p= char(K)$ otherwise.
Note that $a$ is periodic in $L$ iff $a$ is periodic in $L[k]$.

\para{Modular and fixed-field internality}

We return to difference fields.  If we take also derivations the results below regarding the fixed field generate to the constant fixed field of any subset of the derivations; but
we do not know if (I2) holds for modularity in these cases.

Let $\pi_{fix}$ (over $K$) be the set of difference field extensions  $L/K$ generated by 
periodic elements.   Obviously (I1) holds.  In fact (I2') holds too; indeed
$L/K$ is generated by periodic elements iff $L[k] / K[k]$ is.   Let $\pi_{fix-int}$ be the set of
$L/K$ that are qf-internal to $\pi_{fix}$.  It follows easily from the definition that $\pi_{fix-int}$
also satisfies (I1,I2').  

Say $L/K$ is qf-$\pi$-analyzable   if there exists a sequence $K \leq K_1 \leq \ldots \leq K_n=L$
with $K_{i+1}/K_i$ qf-internal to $\pi$.

\<{cor}     
  \<{enumerate}
\item    Let $L \in \mC$, $tr. deg._K L \leq n$.   Let $\Uu$ be a universal domain
containing $L$.  Then $L/K$ is not modular in $\Uu$ iff there exist $K \leq K' <L' \leq L$ with 
$L'/K'$ non-algebraic and in $\pi_{fix-int}$.   
\item  Let $\pi_{mod;K;n}$ be the set of  difference field extensions  $L/K$ 
that are modular in some $\Uu$, of transcendence degree $\leq n$.  Then any $L \in \pi_{mod;K;n}$ is modular in any $\Uu$ containing $L$; and (I1,I2') hold.
\>{enumerate}
\>{cor}
\prf  
Note   that (1) implies (2) immediately.  One direction of (1) is clear: 
If $L/K$ is modular in some universal domain $\Uu$, then ceratinly $K',L'$ cannot exist.
We prove the converse part of (1)  by induction on $n$.  By (2), properties
I1,I2' hold for r $\pi_{mod;K;m}$ for $m<n$.  
  Now suppose $L/K$ is not modular in $\Uu$.  Then by \cite{CH},  $L/K$ is not almost -orthogonal to some $Fix(\tau)$, or else to a 
    modular type $q$ over $K$ of $SU$-rank one.   In the first case by Theorem \ref{6.3}
    there exists $K_1$ with $K_1/K$ non-algebraic and qf-internal to $\pi_{fix}$.  Assume the
    second possibility.    If $q$ has transcendence degree $n$ then $L/K$ is itself of SU-rank one, and modular,
a contradiction.  Otherwise, let $L'$ be the difference field generated by a realization of $q_{qf}$.  Then $L' \in \pi_{mod;K';<n}$, and $L/K$ is not orthogonal to $L'$. Since (I1,I2) hold for
  $\pi_{mod;K';<n}$ there exists $K'_0 \leq L$, $K \neq K'_0$, with $K'_0/K$ qf-internal to
  $\pi_{mod;K';<n}$.  So $L/K'_0$ is not modular.  Again using induction we may find 
  $K'_0 \leq K' < L' \leq L$ as required.   \eprf  
  
  %Let $\pi_{mod;K} = \union_n \pi_{mod;K;n}$.

\<{cor}   Let $L \in \mC$.  Then there exists a finite sequence
$K=L_0 \leq L_1 \leq \ldots \leq L_n=L$ such that for each $i$,   $L_{i+1}/L_i$ is $L_i$-primitive and either finite, or  qf-internal to $\pi_{fix}$, or else it is modular and internal to   $\pi_{mod;K} = \union_n \pi_{mod;K;n}$.  \>{cor}

\prf   Immediate from the fact that $\pi_{fix}$, $\pi_{mod;K;n}$ satisfy I1,I2.   \eprf

\section{Limited sets}     \vlabel{limited}

Let $k$  be a  field, $L$ a finitely generated extension field,
and $V$ a variety over $L$.   The original  formulation of the problem we solve requires the notion of
a {\em limited }  subset of $V(L)$.  
Geometers are familiar with constructions interpreting varieties 
over fields such as $k((t))$, with inductive and projective systems of varieties over $k$.
See for instance \cite{[BL]}, \cite{[GK]}.   In this language, viewing $V(L)$ as an Ind-variety over
$k$, a subset of $V(L)$ is {\em limited} if it is contained in a finite-dimensional $k$-subvariety of $V(L)$.   We now give an account of this for model theorists; see \cite{kam} for a detailed 
explanation of these ideas.  

Let $k$ be a structure.  We will say {\em constructible} for quantifier-free
definable over $k$.   

A structure $N$ for a finite relational language
$L$ is  {\em piecewise-definable} over another structure $k$ if there exist constructible $L$-structures
$N_i$ over $k$ and definable $L$-embeddings $N_i \to N_{i+1}$ such that $\lim_i N_i$ is isomorphic to $N$.   A {\em $k$-definable subset} of $\lim_i N_i$ is just a  definable subset of some $N_i$.   Similarly for {\em piecewise-constructible}.

\<{lem} \vlabel{limited-0} Let $N = \lim_i N_i$ be piecewise-constructible over $k$, 
and let $S$ be a quantifier-free definable subset of $N^k$.  Let $\alpha_i: N_i \to N$ be
the canonical map, and write $\alpha_i$ also for the induced map $N_i^k \to N^k$.  
Then $\alpha_i \inv (S)$ is a constructible subset of $N_i^k$.  

Let $N'$ be an $L$-structure,   quantifier-free definable over $N$.  If    $N$ is piecewise-constructible over $k$, then so is $N'$.
\>{lem}

\prf  Immediate from the definitions.  \eprf

We now consider fields and $k$-algebras.  We take $L_k$ to be the relational language
with a relation for any $k$-constructible set.

\<{lem} \vlabel{limited-1}

 (1)  Let $k$ be a field, and let $L=k(b_1,\ldots,b_n)$ be a finitely generated field extension of $k$.  
Then $(L,+,\cdot,b_1,\ldots,b_n,k)$ is piecewise constructible over $k$ (with parameters in $k$.)    More precisely there exists a piecewise constructible
$k$-algebra $L'$ and an isomorphism $\psi: L  \to L'$ of $k$-algebras.  

(2)  For any variety $V$ over $L$, $V(L)$ can be viewed as   piecewise-constructible
over $k$.  (I.e. $\psi (V(L))= V^{\psi}(L')$ is piecewise-constructible over $k$.)

 \>{lem}

\prf   (1) is clear.  For (2), $L$ is a finite extension of  a purely transcendental extension $k(t)=k(t_1,\ldots,t_n)$ of $k$. Clearly $L$ is  quantifier-free definable over $k(t)$.
Hence by (1) it suffices to show that  $k(t)$  is piecewise-constructible
over $k$.  Indeed let  $S_n$
be the set of rational functions $f(t)/g(t)$ with $\deg(f),\deg(g) \leq n$, and let 
$+,\cdot $ be the graphs of addition and multiplication restricted to $S_n^3$.  Then 
$\lim_n S_n = k(t)$.   

(3) Note that the $k$-algebra isomorphism $\psi$ induces a map $V(L) \to V^{\psi}(L')$,
also denoted $\psi$.  (3) follows from (1) and (2).  \eprf

  Let $L'$ be constructed as in 
  \ref{limited-1} (1); we can view $L'$ as the union of an increasing system of 
  $k$-constructible $L_k$-structures $L'_i$.   Let $b$ 
  be a finite tuple of generators of $L'$ over $k$.   For any $n$, 
  $Y_n(b)=\{f(b)/g(b); \deg f, \deg g \leq n, g(b) \neq 0 \}$ is contained in some $L'_i$.
 Conversely, any $L'_i$ is contained in some $Y_n(b)$. It follows that any $k$-algebra
 automorphism of $L'$ preserves the family of sets contained in some $L'_i$.

{\bf \noindent Definition.}    Let $L$ be a finitely generated extension field of $k$, and let $V$ be 
a variety over $L$.    A subset $Y$ of $V(L)$ is
called {\em limited} if for some isomorphism $\psi: L \to L'$
to a piecewise-constructible $k$-algebra as in \ref{limited-1}, $\psi(Y)$ is contained in an $k$-constructible subset of 
the piecewise-definable set $V^\psi (L')$.

If $L$ is a piecewise-constructible $k$-algebra, $V$ a variety over $L$,  the $k$-constructible subsets of
$V(L)$ can be characterized as the images under piecewise rational constructible functions defined over $L$, of a constructible set over $k$.  This class is preserved
under automorphisms of $k$-algebras.  In the above definition, one can therefore
replace ``some'' by ``every isomorphism $\psi: L \to L'$.
 
 \<{lem}  \vlabel{limited-char} 
 $Y \subseteq V(L)$ is limited iff there exists a constructible set $U$ over $k$, and a constructible map $g: U_L \to V$, such that  $Y  \subseteq g(U(k)) $.  
 If $Y$ is limited, one can choose $g$ to be injective on $U(k)$.  
 \>{lem}  
 
 \prf  We have an isomorphism  $\psi \inv:  \union L_i \to L$, where $\union_i L_i$ is a  piecewise-definable field over $k$.  This induces $\psi \inv: \union V(L_i) \to V$.  Now $V(L_i)$ 
 is  a constructible set over $k$, and $h_i \psi \inv | V(L_i)$ is a constructible map.  By definition, if
 $Y$ is limited then $Y$ is contained in the image of one of these maps.   This gives
 $g,U$ with $Y  \subseteq g(U(k)) $.  By Lemma \ref{limited-0},  $g \inv (=)$ is a constructible
 equivalence relation on $U(k)$.  Factoring it out,   we may  take $g$ to be injective on $U(k)$.   
 
 \eprf

Let us mention two further equivalent formulations, one geometric and one model-theoretic.  

(1)  Assume $V$ comes with a projective embedding and a notion of   height
applies, cf. \cite{lang}.      Then a limited subset of $V(L)$ is a set of bounded height.      This equivalence is standard.  For instance over $k(x)$, a point of $\Pp^n(k(x))$ can be written
in projective coordinates as $(f_0(x): \ldots : f_n(x))$ with $f_i$ polynomials without common
factors, and then the height is the maximal degree of $f_i$.  More generally see
\cite{lang} Chapter 3, Proposition 3.2, or the box below it; 
and recall   that the set of rational functions on $W$ whose polar divisor is bounded
by some fixed divisor forms a limited set, in fact a finite-dimensional $k$-space.  Compare also
Lemma \ref{heights-limited}.  

(2)  Let $T$ be the $\omega$-stable theory of pairs $(k,K)$ of algebraically closed fields, with $k < K$; cf. \cite{Kei}.   The completions of $T$ are obtained by specifying the characteristic.  
%$(K,k)\subset (L,\ell)$ are models of $T$, then $(K,k)\prec (L,\ell)$ if
%and only if $\ell$ and $K$ are linearly disjoint over $k$. 
If $D\subset
k^N$ is definable (with parameters) in the $\call_Q$-structure $(K,k)$,
then it is $k$-constructible. 

  Assume $(k,K) \models T$.  Then 
 A subset $Y$ of $V(K)$ is {\em limited} if it is $k$-internal,
i.e. $Y \subseteq \dcl(b,k)$ for some finite $b$.   In this case $Y  \subseteq V(L)$
 for some subfield $L$ of $K$,
   finitely generated over $k$, and $Y$ is limited in the sense defined above for 
   finitely generated extensions.

 So far, the discussion involved only algebraic  varieties over fields $k,K$, and no dynamics.
Now assume given also a subvariety $S$ of $V \times V$.  We continue to assume:
$k$ is algebraically closed, $K$ a finitely generated extension field of $k$,
$S,V$ are defined over $K$,

 \<{lem} \vlabel{limited2} The following are equivalent:

 (1)  There exists a possible reducible variety $V'$ over $k$, $S' \leq (V')^2$
 with quasi-finite, dominant projections to $V'$, 
  and a dominant rational map
 $V' \to V$ carrying $S'$ to $S$.   
   
 (2)  For some limited subset $Y$ of $V$, for any $n$ and any nonempty open  subvariety
 $W$ of $V$, there exist $a_0,\ldots,a_n  \in W \meet Y$, with $(a_i,a_{i+1}) \in S$.

 \>{lem}

\prf  Assume (1).   
Let $\phi: V' \to V$ be a dominant rational map taking $S'$ to $S$.   Let $W$ be as in (2).
Let $W' $ be the set of points $v \in V'$ such that $\phi(v)$ is defined and $\phi(v) \in W$.
Choose  $b_0,\ldots,b_n \in W'(k)$, $(b_i,b_{i+1}) \in S'$; and let $a_i = \phi(b_i)$.

Now assume (2).  We may assume here that $k$ is saturated (i.e. has infinite transcendence degree over the prime field.)  In this case by compactness there exist $a_i \in Y$ for $i \in \Zz$,
with $a_i$ avoiding any $K$-definable proper subvariety of $V$, and with $(a_i,a_{i+1}) \in S$.
Let $U,g$ be as in Lemma \ref{limited-char}, with $g$ to be injective on $U(k)$.   
We will use the following general principle:  
 
(*) for any constructible $W \subseteq V^2$, the pullback $g \inv (W) $ is a constructible subset
of $U$.  
By Lemma \ref{limited-0},  
 we obtain a constructible $S'' =g \inv (S) \leq U^2$,
 such that $g$ carries $S'$ to $S$.  The projections $S'' \to U$ have finite fibers
 since this is true for the projections $S \to V$.  Let $a_i' = g \inv (a_i)$.  
 We may take $U$ to be a finite union of varieties,
and $g$ piecewise rational (pre-compose with a Frobenius power if necessary).  Let $V'$ be the Zariski closure of $\{a_i': i \in \Zz\}$.  Let
$S'= S'' \meet (V')^2$.  Then the projections $S' \to V'$ are dominant, since their image
contains all the $a_i'$.

 \eprf

 \<{remark}  \vlabel{limited3} Suppose this weakening of (2) holds:  $S \leq V^2$ has quasi-finite projections;
 and 
 
 (2') For some limited subset $Y$ of $V$, for any $n$, 
    there exist distinct $a_0,\ldots,a_n  \in  Y$, with $(a_i,a_{i+1}) \in S$.   
    
  Equivalently, after saturating we obtain
  
  (2'')     there exist distinct $a_i  \in  Y (i \in \Zz)$, with $(a_i,a_{i+1}) \in S$.   
 
 Let $V'$    be the Zariski closure of $\{a_i: i \in \Zz\}$, and $S' = S \meet (V')^2$.
 Note that $S'$ projects dominantly and quasi-finitely to $V'$; moreover now (2) holds.

We will thus investigate the consequences of (2); assuming (2') instead, they will automatically
apply to  some infinite $(V',S') \leq (V,S)$.
 
 \>{remark}

This can also be reformulated using canonical heights, when they are available.
Assume $S$ is the graph of a morphism $s$.   
Assume   ${h}$ is a function on $V(K)$ such that
 \<{description}
\item(i) A subset $Z$ of $V(K)$ is limited iff ${h}(Z)$ is bounded above in $\Rr$.
\item(ii)  For some $\kappa>1$, 
  if $a \in V(K)$ and $b=s(a)$ then 
   ${h}(b) = \kappa {h}(a)$.   
\>{description}
     Lemma \ref{limited3} (2') is equivalent to the statement:
 for some $\e>0$, for any $n$, there exist 
  distinct $a_0,\ldots,a_n  \in  V(K)$, with $(a_i,a_{i+1}) \in S$, and ${h}(a_i) < \e$.
If this holds, then ${h}(a_0) < \e \kappa^{-n}$.  Conversely if just $a_0$
can be found with ${h}(a_0) < \kappa^{-n}$, letting $a_i = s^i(a_0)$,
we have ${h}(a_i) = \kappa^{-i} \leq 1$.
% so $h_W(a_i)  <  ||h-h__W||_\infty + 1$.  Hence
(2'') is   equivalent in this situation to:  
 
(2h)  For any $\e>0$ and $n$,  there exist   $a_0 \in V(K)$ with $s^n(a_0) \neq 0$ 
and  ${h}(a_0) < \e$.
 
If one knows that for each $n$ there are only finitely many fixed points of $s^n$,  the condition $h^n(a_0) \neq 0$ can be replaced by:

 There exist infinitely many  $a_0 \in V(K)$ with   ${h}(a_0) < \e$.  
 
 This is the formulation used in \cite{baker}.

 \proof of  \ref{g2}.     The finite version follows by compactness from the
 qualitative one.  So assume some $\phi$-orbit $a_1,a_2, \ldots$ is contained
 in a limited subset of $V(L)$.   View $V(L)$ as a direct limit of constructible sets  $V_i$ over $K$.
 $\phi$ induces a function $F: V'(L) \to V(L)$, ($V'$ being the domain of definition of $\phi$);
 this is a morphism of Ind-constructible sets; i.e. for each $i$, for some $j$,
 $F|V_i$ is a constructible function $V_i \to V_j$.    By assumption  
 $\{a_1,a_2,\ldots \} 
 \subseteq  V_i(K)$ for some $i$.  Let $U'$ be the Zariski closure in $V_i$
 of this set.  Then $F(U') \subseteq U'$.  For some $m$, and some component $U$ of $U'$ of maximal dimension, we have $F^m(U) =U$.  We can view $\bU=(U,F^m) \in AD_K$.
Then $\bU_L$ dominates $(V,\phi^m)$.  By Theorem 3.3 of Part II of this paper, $(V,\phi)$ descends to $K$. \eprf

\<{rem}   The proof goes through for difference varieties, not necessarily algebraic dynamics.
If the difference variety is generated by a relation $(a,b) \in S$, with $S$ a correspondence,
we conclude that there for any limited subset $Y$ of $L$, for some $n$  are no $a_1,\ldots,a_n \in L$ with $(a_i,a_{i+1}) \in S$.  \>{rem}  

\<{question}  if 
$(V,\phi)$ is a field-free algebraic dynamics over $\Qq$, or over $K(t)$, do there
exist $a,b \in \Qq^{alg}$ or $K(t)^{alg}$ in the same orbit of Galois, and with $\phi(a)=b$? \>{question}

It would follow that $a$ is periodic.  When  $(V,\phi)$ is polarizable, Fakhruddin has shown that
at least periodic points exist (and are Zariski dense.)

 \def\uV{\underline{V}}  \def\uU{\underline{U}} \def\uW{\underline{W}} \def\uR{\underline{R}}
    
\sect{Canonical heights}\vlabel{heights}

 Perhaps the most characteristic feature of algebraic dynamics is the {\em canonical height} associated to an algebraic dynamics  $(X,\phi)$ over a field $K$; see
  \cite{lang}, \cite{HS},   \cite{CS}.    The descent questions treated in the present paper are usually stated in terms of canonical height.  We review this concept here in the case
 of a a function field $K=k(C)$, with $C$ a curve over $k$.   Our presentation is more general than the usual one in three ways, compatible with the more general setting of the paper: (1) we allow $k$ to be a difference field, so that $\phi:X \to X^\si$; this generalization requires no effort, and the reader interested in the basic case may take
 $\si=Id$.  (2)  we allow ``probabilistic" dynamics given by correspondences, rather than morphisms  or rational maps; (3)   in the classical case $\si=Id, \phi$ a morphism, we give a construction that does not depend on a polarization assumption; and show
that it is nontrivial under much weaker conditions than polarizability.

The dynamics will take place on   complete normal  varieties $V$ over $K$.  
Such a variety can be viewed as the generic fiber of a complete normal variety $\bV$ over $C$,
i.e. a variety over $k$ with a dominant morphism $j:\bV \to C$; 
we fix $\bV$ too.  In this section we use boldface for varieties over $C$; difference
varieties will be considered only in the explicit geometric form $(V,\phi)$ or $(V,S)$.

   The basic properties of 1-cycles and divisors are reviewed below.   For any variety $V$ let  $NS(V)=Pic(V)/Pic^0(V)$ denotes
   the group of Cartier divisors on $V$, up to algebraic equivalence; see \cite{langneron}.
   
We assume $V$ and $\bV$  are smooth.   It follows that the notions of Weil and Cartier divisors
are the same on $V$ or $\bV$; we will refer to them as divisors.  The coincidence of Weil and Cartier  is used to
show that the natural map $NS(\bV) \to NS(V)$ is surjective, and to define a pushforward
$NS(S) \to NS(V)$ when $S \to V$ is finite.   Together with a certain
statement on Galois covers in case the dynamics is not rational, this will be our only use
of the smoothness assumption.  

We assume $k$ and $K$ come with compatible endomorphisms $\si$, and let $V'=V^\si$;   the reader may take $\si=Id$ if desired.  Let $S \leq V \times V'$ be a subvariety, with
finite projection $p: S \to V$ of degree $d$, and generically finite projection $q: S \to V'$ of degree $d'$.  We make no smoothness assumption
on the dynamics $S$.  We will explain
below how an action $S_* = \si \inv q_*p^*$ is induced on 0-cycles over $V$. 
  We can
view $S_*$ as a non-deterministic dynamics on $V(K^{alg})$, going from $a$ to $b_i$ with probability  $m_i/d$ if $S_*(a) = \sum_{i=1}^d m_i b_i$.  
 In case
$S$ is the graph of a morphism $\phi$, this is compatible with the action of $\phi$
on $V(K^{alg})$.   Similarly,  If $h$ is a function 
on $V(K^{alg})$, and $c=\sum m_i a_i$ is a 0-cycle, we let $h(c)=\sum m_i h(a_i)$.

\<{defn}  \lbl{ch} A {\em canonical height} is a  a finite dimensional $\Rr$-vector space $\Lambda$, an expanding linear transformation $\lambda: \Lambda \to \Lambda$  and  a   function $h: V(K^{alg}) \to  \Lambda $   such that 
  for  any $0$-cycle $a$ on $V$,
  $$h(S_*(a)) =   \lambda (h (a))$$ \>{defn}  
  
%If $v \in \Lambda$ is an eigenvector of $S_*$ with real eigenvalue $\lambda_v$, 
%we can write  $h(a)(v)= h_a v$; then we have $\lambda_v>0$, and $h_a$ is a canonical
%height in the sense of Call and Silverman, \cite{CS} (when $S$ is the graph of a morphism.)
%  

By {\em expanding} we mean that every complex eigenvalue lies outside the unit circle.
Classically (cf. \cite{CS}, discussed below), one only takes the case $\dim(\Lambda)=1$.

Classically (cf. \cite{CS}, discussed below), one only takes the case $\dim(\Lambda)=1$.

By complexifying and taking eigenvectors $v_\nu$,  with eigenvalues $\lambda_\nu$ we could  obtain 
maps $h_\nu$ into $\Cc$ with $h_\nu(S_*(a)) =   \lambda_\nu (h_\bu (a)) $ 
and  
take absolute values to get $|\nu|$ into $\Rr$, with $|\lambda_\nu| > 1$.  This may look closer to the classical case.  However the
construction of $\Lambda$ will be more canonical, and gives more information in the non-semi-simple case, so we do not choose eigenvectors.   At all events   we do not assume the existence of real eigenvalues.
 
\<{thm} \lbl{ch1} There exists a canonical height $h: V(K^{alg}) \to \Lambda$, canonically associated to  
$(V,S)$. \>{thm}

We think $\Lambda$ is rarely trivial.  We formulate two statements in this direction
when  $S$ is the graph of a rational map $\phi$.  One is simply that $\phi$ is not 
birational.  Another, applying even in the birational case,   links with the Northcott results of this paper (e.g. \ref{1.9}).   Assume  $\phi: V \to V^\si$ is a morphism; then $\phi$ induces
a homomorphism $NS(V^\si) \to NS(V)$, and by composing with $\si \inv$
we obtain an endomorphism $\phi^*$ of $NS(V)$.

\<{thm}  \vlabel{thmh}    Let $V$ be a smooth projective variety over $K$, and   $S$ is the graph of a rational map $\phi: V \to V^\si$. 

a)     Assume $\phi$ is not birational.  Then the  canonical height $h$ of Theorem \ref{ch1} is nontrivial.

b)   Let $\phi: V \to V^\si$ be a morphism, and 
 assume no non-torsion element of $NS(V) $
is fixed by a power of $\phi^*$.   Then  $h(a) = 0$ iff  the $\phi$-orbit of $a$ is contained in a limited set. \>{thm}

 In fact if in (b) the morphism $\phi: V \to V'$ is   finite,  then any ample height  is bounded uniformly on the set of all elements of canonical height $0$.
 
Note that if $(V,\phi)$ is primitive,  non-isotrivial, and has no periodic subvarieties of positive dimension, then by \ref{g2} the  condition in (b)
is equivalent to pre-periodicity of $a$.  Conversely, if $a$ is pre-periodic then 
any canonical height must vanish on it.  Hence in this case at least $h$ is the universal
canonical height, i.e. any canonical height $h' : V \to \Lambda'$ is the composition of $h$ with    a linear map $\Lambda \to \Lambda'$.  

%By Corollary \ref{1.9}, if $\phi$ does not constructibly descend to $k$ and is field-free,
%then the $\phi$-orbit of such an $a$ cannot be Zariski dense.
%  

This construction of canonical heights would formally extend to  number fields $K$ given
some Arakelov theoretic information, of which we are uncertain.  We 
would need an Arakelov variety $\bV$ with $\bV_K=V$, a definition of 
$NS(\bV)$ (with real coefficients) and a pairing $A_0(\bV) \times NS(\bV) \to \Rr$, and
a surjective $\Rr$-linear (``Gysin'') map $NS(\bV) \to NS(V) \tensor \Rr$.   Moreover
given a correspondence $S \leq V \times V^\si$, we require pullback maps $NS(\bV) \to NS(\bS)$, and push-forward maps for 1-cycles, related by a projection formula.   No finiteness statement on $NS(\bV)$ is needed.
 
We begin with a review of (non-dynamical) heights over function fields.
%   Over a function field, there are two approaches.
%We may view the   field geometrically as the function field $K=k(C)$ of a curve (or more general variety) $C$, and the given variety $V$ as the generic fiber of a variety $\bV$
% over $k$, with a morphism $\bV \to C$.  Then a point of $k(C)^{alg}$ is seen as a curve
% on $\bV$, and heights are connected to intersections with divisors.  The second
%  approach (N\'eron) begins with a set of absolute values of $K$, satisfying the product formula.  The relation between the two is not tautological.    The geometric approach has  the advantage of giving a clear visual image, and it is this path that we choose here.  

% A technical note: we will assume $V$ is smooth, and in fact is the $K$-fiber of a smooth complete variety $\bV$ over a   $k$.   
%  We will actually only use intersections with divisors; intersections with Cartier divisors require no smoothness.
% However the notion of ``divisors" that we use needs to be closed under proper pushforwards,
% hence we use Weil divisors.   It would suffice to assume, in place of smoothness,  that the notions of Weil and Cartier divisors coincide on $V$, perhaps even only up to algebraic equivalence.   At all events, smoothness of the ambient algebraic variety is not an overly restrictive assumption; in characteristic $0$ it can be achieved with at worst a birational change to the variety; we impose no smoothness condition on  the dynamics.  

\para{\bf   Cycles.}  
  Fix a   complete  variety $W$ of dimension $n$ over an algebraically closed field.  
 We outline the most basic concepts of intersection theory on $W$; for this purpose
 we take $W$ to be smooth and discuss all dimensions, but we will really use only cycles of
 dimensions 0,1 and $n-1$.   See a similar summary in \cite{ha}, or   \cite{fu}.     

Let $U,V$ be subvarieties of $W$ of complementary dimension.
If the intersection $U \meet V$ is transversal, we write $U \cdot V$ for  the number of intersection points.   From this basic geometric data one forms an algebraic structure as follows:

Let $C_l(W)$ be the Abelian group freely generated by the
$l$-dimensional irreducible subvarieties of $W$.
 if $U = \sum_{i=1}^m \a_i U_i \in C_l(W)$, 
$V= \sum_{i=1}^m \beta_i V_i \in C_{n-l}(W)$  
 and each intersection $U_i \meet V_j$ is transversal,
we say that $U,V$ are transversal and define $U \cdot V = \sum_{i,j} \a_i \beta_j U_i \cdot V_j$.  
This symbol is then extended to the non-transversal  case, as follows.  Two 
cycles $U,U' \in C_l(W)$    are said to be  {\em numerically equivalent}  if for any  $V \in C_{n-l}(W)$ transversal to $U$ and to $U'$, we have $U \cdot V = U' \cdot V$.  Write $[U]$ for the class of $U$ 
up to numerical equivalence, and let  $A_l(W)= \{[U]: U \in C_l(W) \}$ be the quotient of  
$C_l(W)$ by the cycles numerically equivalent to $0$.  
Write $A^l(W)$ for $A_{n-l}(W)$.  For reducible subvarieties $U$ of dimension $l$, we let $[U]=\sum_{V} [V]$,
where $V$ ranges over the irreducible components of $U$ of dimension $l$.
There exists a unique 
bilinear pairing $A_l(W) \times A^{l}(W) \to $ such that  $[U] \cdot [V] = U \cdot V$ when the right hand side is transversal.

 If $l_1+l_2+l_3=n$, it  is similarly possible
 to define a trilinear pairing $A_{l_1} \times A_{l_2} \times A_{l_3} \to \Zz$, 
 $(U_1,U_2,U_3) \mapsto U_1 \cdot U_2 \cdot U_3$, agreeing with the number of intersection points 
 in the transversal case.  In fact there exists a bilinear map $A_{l_1} \times A_{l_2} \to A_{n - l_1-l_2}$ such that $ U_1 \cdot U_2 \cdot U_3 = (U_1 \cdot U_2) \cdot U_3$.

Let $W'$ be another smooth complete variety, and  $f: W \to W'$  a   morphism.   For an irreducible  subvariety $U$ of $W$, define
  $f_*([U]) = \deg(f|U) [U'] \in C_l(W')$, where 
$U'=f(U)$ the image in $W'$, and $\deg(f|U)$ is defined to be field extension degree
$[k(U):k(U')]$ if this is finite, $0$ otherwise.  Extend by linearity to a
linear map $f_*: C_l(W) \to C_{l}(W')$.  Then $f_*$ induces a homomorphism $f_*: A_l(W) \to A_l(W')$.  

$A_0(W)$ can be identfied with $\Zz$ via the degree map  $\sum \a_i [u_i] \mapsto \sum \a_i$.  The push-forward $f_*$   preserves degree on $A_0$.

When $W = \Pp^n$ and $V$ is a hyperplane, $U \cdot V$ is the projective degree of $U$.
For any integer $m$, the family of all curves $U$ on   $W$ with $U \cdot V \leq m$ is therefore a limited family.

For our purposes a {\em divisor} is an element of $A^1(W)$.   By N\'eron-Severi, cf. \cite{langneron}, this is a finitely generated Abelian group.   We will also write $NS(W)$ for $A^1(W)$, and $NS_\Rr(W)$
for $\Rr \tensor NS(W)$.  $NS$ is a contravariant functor for surjective morphisms. 
We have the  projection formula (\cite{fu}, 2.3(c))):   given $f: W \to W'$, 
 for $U \in A_l(W), D \in NS(W')$ we have: $f_*(U \cdot f^*(D)) = f_*(U) \cdot D$.   
In case $l=1$, using our identification of $A_0$ with $\Zz$, and the fact that $f_*$
preserves degrees, this also reads:  $U \cdot f^*(D) = f_*(U) \cdot D$.   

A divisor $D$ on $W$ is called {\em very ample} if it is the pullback of a hyperplane, under
some projective embedding of $W$; {\em ample}   if $mD$ is very ample for some $m > 0$.  Such a divisor $D$  inherits the property noted for the hyperplane divisor on $\Pp^n$:  the family
of curves $U$ on $W$ with $U \cdot D \leq m$ is a limited family.

 We can define heights using {\em either} intersections of subvarieties $S$ on $W=V \times V'$
 with certain divisors and curves   from $V$ and $V'$, 
{\em or} using intersections of divisors and curves on $S$ itself.   The latter
is more efficient since it can be defined for Cartier divisors without assuming $V$ is smooth.  However,
we need to push forward Cartier divisors under generically finite morphisms, and cannot
do it unless Cartier divisors coincide with Weil divisors, at least up to algebraic equivalence.  So at all events some smoothness
assumption is needed.    At all events, smoothness of the ambient algebraic variety is not an overly restrictive assumption; in characteristic $0$ in particular,  it can be achieved with at worst a birational change to the variety, and then $\bV$ can be chosen smooth once $V$ is.
   We impose no smoothness condition on  the dynamics $S$.

\para {\bf Weil Height}.  \vlabel{weilh}

We   work with the   data $k,K=k(C), j: \bV \to C$  with generic fiber $V$, as above.
 Cycles on $\bV$ will always be assumed to be defined   over $\k$.  
 If $U \leq \bV$ is
an $l$-dimensional irreducible variety defined over $\k$, then $U \meet V$
is either empty or an $l-1$-dimensional $K$- irreducible subvariety  of $V$.
Any $K$-irreducible subvariety of $V$ can be written uniquely in this way.   This gives 
  a homomorphism
$$\rho:  C_l(\bV) \to C_{l-1}(V)$$
whose kernel is generated by the classes of subvarieties that project onto a non-Zariski
dense subset of $C$ (i.e. a finite subset of $C$.)

In the opposite direction, given $A \in C_{l-1}(V)$ there exists a unique $\bA \in C_l(\bV)$
whose support has no component projecting to a point of $C$, and with $\rho(\bA)=A$.
We denote $\beta(A)=\bA$.  
  
Up to Galois conjugacy, a point   of $V(K^{alg})$ can be identified with an
 irreducible element $a$ of $C_0(V)$, i.e.
a cycle with non-negative coefficients, not all zero, which is not the sum of two other such.
In this case $\bba=\beta(a)$ is a curve on $\bV$, the morphism  $j|\bba: \bba \to C$ is finite, and has
degree equal to the degree $d(a)$ of $a$ as a 0-cycle.

Let $D$ be a divisor on $V$, and let 
 $\bD$ be a divisor on $\bV$, defined over $k$,  restricting 
to $D$ on $V$.   
 
Let $a \in C_0(V)$. 
We define the {\em Weil height } of $a$ by:  $  h_{ \bD}(a) = \beta(a) \cdot \bD / d(a) $.
Here $\beta(a) \cdot \bD$ is the intersection number.    Note that
$h_{\bD}(na)= h_{\bD}(a)$, so $h$ factors through   the projectivization
$\Qq \tensor C_0(V) / \Qq^*$.  

 If $\bV  = C \times \Pp^n$ and $D$ is the hyperplane divisor on $\Pp^n$ pulled back to $\bV$, then we have the usual Weil height over $k(C)$ on $\Pp^n$.  
%% 

%%  The definition is more familiar in local terms.   See \cite{lang}.  
%%  To make the connection
%%consider the case:  $V=\Pp^1, D=\{\infty\}=[(1:0)],  \bV = C \times \Pp^1$, $(a_0:a_1) \in \Pp^1(C) \setminus \{\infty\}$.  Then $a=a_0/a_1$ is  a rational function on $C$, and $\beta(a)$ is the Zariski closure of
%%the graph of $a$.  We have $\bD = C \times \{\infty\}$.  Now the intersection $\beta(a) \cdot \bD$
%%can be expressed as a sum of the actual intersection points $p_1,\ldots,p_k$ of $\beta(a) \meet \bD$, 
%%with certain multiplicites.  Let $v_p$ be the valuation on $K=k(C)$ corresponding to the point $p$.
%%Then it is clear that $(p,\infty) \in \beta(a) \meet \bD$ iff $v_p(a) <0$.  It can be shown that $-v_p(a)$
%%is in fact the intersection multiplicity.  Thus $h_{ \bD}(a) =\sum_{p} \max{(-v_p(a) ,0)} $.
%%Using the product formula $\sum v_p(a_1)=0$, this can also be written:  $  h_{ \bD}(a) =\ \sum_p \max(v_p(a_0),v_p(a_1))$.
%% 
      
\<{lem}   \vlabel{height1}   If $\bD'$ is another divisor on $\bV$, restricting to $D$ on $V$, then 
$h_{ \bD} - h_{\bD '}$ is a bounded function on $V(K^{alg})$.  \>{lem}

\prf  Let $L = \bD - \bD'$; write $L=\sum m_i L_i$ where $L_i$ is an irreducible hypersurface of $\bV$; since $\bD,\bD'$ agree on a generic fiber of $j$, we have $j(L_i)$ non-Zariski-dense
for each $i$, i.e. $j(L_i)$ is finite.  So $L_i$ is contained in a fiber of $j$.  
Let $a \in V(K^{alg})$, $\bba = \beta(a)$.  
For a generic
fiber, hence for each fiber $L'$ of $j$ we have $\bba \cdot L' = d(a)$.  
  It follows   that  $ | \bba \cdot L |  \leq  \sum |m_i| d(a)$.    So $|h_{\bD}(a) - h_{\bD'} (a)| \leq \sum |m_i|$.  \eprf  \smallskip

Let ${\mathcal F}(V)$ be the space of functions $V(K^{alg}) \to \Rr$, modulo the bounded functions.
The 
  class   $h_{\bD}$ in ${\mathcal F}(V)$ depends only on $D$; we denote it $h_D$.
%We thus write $h_D$ for   $h_{\bD}$ whenever we refer to a property of $h_{\bD}$ that
%is invariant upon adding a bounded function.  

\smallskip
\para{\bf \noindent Bounded height and limited families}   We recall (cf. \cite{langneron}, Property 1F)):  

{\bf \noindent Lemma} \vlabel{heights-limited} Assume $D$ is very ample.  Fix $d_0 \in \Nn$, and $\a \in \Rr$.  Then 
$$\{a \in V(K^{alg}): d(a) \leq d_0,  h_{\bD}(a)  \leq \a \}$$ is a limited set.

\prf   Using Lemma \ref{height1}, we may replace $\bD$ by any divisor on $\bV$ restricting to 
$D$ at the generic fiber.  So we can add to $\bD$ any divisor whose projection to $C$ is not Zariski dense. 
 The linear system $L(D)$ of rational functions on $V$ with poles at most at $D$
contains elements $f_0,\ldots,f_l$, 
such that $v \mapsto (f_0(v): \ldots: f_l(v))$ is an embedding $V \to \Pp^l$.  These $f_i$
can be taken to be defined over $K$; they are restrictions of functions $F_i$ on $\bV$.  By adding to $\bD$ some divisors whose projection to $C$ is finite, we may assume
$F_i \in L(\bD)$.  By further adding  to $\bD$ a divisor $j^*(D_C)$, where
$D_C$ is a very ample divisor on $C$, we may assume there exist functions $G_i = g_i \circ j$
in $L(\bD)$, such that $c \mapsto (g_0(c): \ldots: g_{l'}(c))$ is a projective embedding of $C$.
In particular $G_0,\ldots,G_{l'}$ do not simultaneously vanish on $\bV$.  
Let $J(v) = (F_0(v): \ldots: F_l(v): G_0(v): \ldots : G_{l'}(v))$.   Then $J$ is a morphism
$\bV \to \Pp^{l+l'+1}$.  If $J(v)=J(v')$ then $j(v)=j(v')$; and for some proper subvariety
$W$ of $C$, if $j(v) \notin W$, then $J(v)=J(v') $ implies $v=v'$.  Now given $a \in V(K^{alg})$, $\bba = \beta(a)$,
$J(\bba)$ is a curve in $\Pp^{l+l'+1}$ of projective degree $d(a) h_{\bD}(a) \leq d_0 \a$.
Thus the $J(\bba)$ span a limited family.  It is clear that $\bba$ is determined by 
$J(\bba)$ in a uniformly definable fashion (away from $W$ we have $\bba = J \inv (J(\bba))$.)
So the family in question is also limited.  \eprf

\para{\bf \noindent Effect of $Aut(K^{alg})$ on heights}

  Let $\g: C' \to C$ be a finite morphism of curves, 
$K' = k(C')$. 
and assume there exists a smooth $\bV', j': \bV' \to C'$,
and a morphism $e: \bV' \to \bV$  with $j e= \g j'$.   
 Let $V'$ be the generic fiber of $\bV'$, and $D' = (e|V')^* (D)$.  
Then
$k(C')^{alg} = k(C)^{alg}$, so we can compare the height of  $a \in V(K^{alg})$ from
the point of view of $C$ and of $C'$.  The point $a$ corresponds to a curve 
$F'_a$ on $\bV'$, with $e(F'_a) = \beta(a)$.   Thus both $h_D(a) :=h_D(\beta(a))$
and $h_{D'}(a) : = h_{D'} (F'_a)$ are defined (up to bounded functions.)

\<{lem}  \vlabel{curvechange} Let $\g: C' \to C$ be a finite morphism of curves, and let
$h_D,h_{D'}$ be as above.  Then in the space ${\mathcal F}(V)$ we have:
$$\deg(e) h_D  =  \deg(\g)  h_{D'}$$ 
 \>{lem}

\prf   Choose a divisor $\bD'$ restricting to $D'$ at $V'$.   
So $h_{\bD'} = h_{D'}$
up to bounded functions.  Let $d'(a)$ be the degree of $j' | F'$.
Fix $a$ and let $F'=F'_a,F=\beta(a)$.  Then $e_* [F'] =  \deg(e | F')  F$
By the projection formula (\cite{fu}, 2.3(c)), and since $e_*$ preserves the degree 
of 0-cycles, we have 
$$[F'] \cdot \bD'  = e_*([F'] \cdot \bD') =  \deg(e | F')  F \cdot D  =  \deg(e | F') [F] \cdot [\bD]$$
So $d'(a) h_{\bD'}(a) =   \deg(e | F') d(a) h_{\bD}( a)$.  Now since $j e= \g j'$ we have
$\deg(\g) \deg(j'|F') = \deg(e|F') \deg(j|F)$ and the claim follows.  \eprf
  
For any finitely generated extension field  $K$ of $k$ we have
 the modular function $\delta: Aut( K^{alg} / k ) \to \Qq$.  It can be defined by
 the ratio of field degree extensions:
$\delta(\si) = \frac{[KK^\si:K^\si] }{[KK^\si:K]}$, where
$K^\si=\si(K)$ and $KK^\si$ is the field compositum.

 Assume now that $V$ is defined over $k$, i.e. $\bV = V \times C$.  
Any $\si \in Aut( K^{alg} / k )$ 
  induces a function $\si_V:  V(K^{alg}) \to V(K^{alg})$, and  composition induces
  an action of $\si$ on ${\mathcal F}(V)$.   Lemma  \ref{curvechange}  implies:
 
\<{cor}  \vlabel{eff}  Let $\si \in Aut( K^{alg} / k )$.  Then in ${\mathcal F}(V)$ we have:
  $$h_{D} \circ \si_V = \delta(\si) h_{D}$$
 \>{cor}

%\para{Canonical heights of solutions of difference equations. }     

The same relation thus holds for the canonical height $h$.  
 Corollary \ref{eff} implies in a   variety of cases that any algebraic solution
 of a difference equation $(x,\si(x)) \in S$ has canonical height $0$.   We have
 $h (\si(a)) = \delta(\si) h(a) $ while $h(s(a)) = \lambda h(a)$, so if
  If $\si(a)=s(a)$ and 
  $h(a) \neq 0$, then $h(a)$ 
 then   must be an eigenvector for $\lambda$  with eigenvalue  $\delta(\si)$.  Thus:  
 
\<{cor}  \vlabel{height0}  Assume $V$ is defined over $k$, $a \in V(K^{alg})$ and $k(a) \cong_k k(s(a))$ by
 an isomorphism  $\si$ taking $a$ to $s(a)$.  Then any of the conditions below implies
 that $a$ has canonical height $0$:  
 
 1)    $S$ is the graph of a rational function, and  $\lambda$ has no rational eigenvalues.
  
 2)   $S$ is the graph of a rational function, and $\delta(\si)$ is smaller than any
 real eigenvalue of $\lambda$. 

3)    $\delta(\si) > \spec(\lambda)$.

 \>{cor}

\para{\bf Dynamics of a correspondence}
 
Let $V$ be a  variety defined over $K$, as above.  
 If $K=(K,\si)$ carries a  difference
field structure, leaving $k$ invariant,  let $V' = V^\si$.    (If one wishes to think of $K$ as a field, let $\si=Id$.)
 Let $S$ be a complete variety of dimension $n=\dim(V)$, defined over $K$, and let 
 $p: S \to V$ and $q: S \to V'$ be morphisms.     We assume $p$ is finite, of degree $d$, and that  $q$ is generically finite, of degree $d'$.  For simplicity, we assume $p$ is separable.

We have $p^*:  NS(V) \to NS(S) $ and $q^*: NS(V') \to NS(S)$.  Since $p$ is finite,
we also have $p_*: NS(S) \subset A^1(S) \to A^1(V) = NS(V)$.  We obtain an 
endomorphism $S^t$ of $NS(V)$, namely $S^t (D) =  p_* q^* D^\si$.  Let $S^* = d \inv S^t$.

Similarly we have $p^*: C_0(V) \to C_0(S)$, using finiteness of $p$.  And we have $q_*: C_0(S) \to C_0(V')$.  We obtain an endomorphism $S_*$ of $C_0(V)$, namely $S_*(a) =   q_* p^*(a) ^{\si \inv}$.  

We will consider two cases:  subvarieties of  $V \times V'$ with the projection maps to $V$ and $V'$; or normalizations $S$ of such varieties.  We note that an arbitrary correspondence $S$ gives equivalent dynamics to one of this form.  Let $S,p,q$ be as above.  
  Let $S'$ be the image of $S$ under $(p,q): S \to V \times V'$.  We can define a dynamics
 using $S'$ with the projection maps to $V,V'$.  Now  the map $\pi: S \to S'$ is finite, and
 we have $\pi_* \pi^* = \deg(\pi)$.  It follows that the dynamics given by $S$ and $S'$
 on $NS(V)$ and on $PA_0(V)$ are the same (on $A_0(V)$
 they differ   by a constant multiple $\deg(\pi)$.)   
Hence we can work with subvarieties of $V \times V'$.  

However, it is convenient to have $S'$ normal and Galois over $V$.  Hence we show
how to replace a given   $S' \leq V \times V'$ by $(S,p,q)$ with these properties, and
with $(p,q)(S)=S'$.    Let $\pi: S \to S'$ be the 
normalization of $S'$ in the Galois hull $L$ of $K(S')/K(V)$.  Then $H=Gal(L/K(V))$ acts
on $S$, over $V$.   For $D \in NS(S)$, at least  if $D$ is represented by a divisor with support
not contained in the ramification divisor $Ram(p)$  of $p$, we  have  $p^* p_* D = \sum_{h \in H} h ^* D$.    
%If $D$ has the form $q^* D'$ for $D'$ a Cartier divisor on $K(V')$, we 
%will see using smoothness of $V$ that this hypothesis holds. 
%may   take $D'$ 
%with support not contained in $q_*(Ram(p))$, and then the hypothesis holds.  

When $S' \leq V \times V'$, 
let $\bS'=\beta(S')$ be the unique irreducible subvariety of  $ \bV \times _C \bV$ with $\bS' \meet (V \times V)= S'$.   
So $\dim(\bS')=n+1$.  In general, let $(\bS,p: \bS \to \bV, q: \bS \to \bV^\si$  be a triple with generic fiber $(S,p,q)$.    If $S$ is the normalization of $S'$ as above, we can 
let $\bS$ be a similar normalization of $\bS'$, so that $H$ acts on $\bS$,  extending the action on $S$.

 Let $\bD$ be a divisor on $\bS$.  

 We have an intersection pairing $C_1(\bV) \times NS(\bV) \to \Zz$, and 
also $C_1(\bS) \times NS(\bS) \to \Zz$.   They are compatible via 
\beq{projection} {  p_* x \cdot y = x \cdot p^* y }   \eeq
for $x \in C_1(\bS), y \in NS(\bV)$; and similarly for $q$.

Recall $\bt_V : C_0(V) \to C_1(\bV)$,  $\bt_S: C_0(S) \to C_1(\bS)$.    Let $\bt_{V'}$ be the $\si$-conjugate.

The  ``Gysin" homomorphism $\rho: NS_\Qq(\bV) \to NS_\Qq(V)$ is surjective, using the fact
that every Weil divisor on $\bV$ is Cartier (or, has a multiple numerically equivalent to one.)
Let  $\g_V: NS_\Qq(V) \to NS_\Qq(\bV)$ be a section of this linear map, and similarly
$\g_{V'}$.  We omit the subscripts when possible.

We have 
\beq{H2.1}        \bt_{V'} q_* = q_* \bt_S    : C_0(S) \to C_1(\bV') 
 \ \ ; \ \ \  \  \bt_V p_* = p_* \bt_S: C_0(S) \to C_1(\bV)      
     \eeq
     
  Note that on $\bS$ the morphism $p$ may not be finite, and so $p^*, \beta$
  need not commute.

Given two functions $f,g$ of a variable $x$ in $C_0(V)$ and possibly other variables $y$,
we write $f \sim g$ if for any $y$, $f(x,y) - g(x,y)$ is bounded on $C_0(V)$.
For elements $b,b' \in NS(\bV)$ we write $b \sim b'$ if $h_b \sim h_{b'}$.

Let $\xi \in C_0(S), y \in NS(V')$.  From \eqref{projection} and  \eqref{H2.1}
 
 we obtain: 
\beq{H3}{    \beta(\xi) \cdot q^* \g (y)    \sim  \beta(q_* \xi) \cdot \g (y)  }   \eeq

\<{lem}  \lbl{transpose1} Let  $x \in C_0(V)$, $y \in NS(V')$.  
Then $\beta(q_* p^* x) \cdot \g (y) \sim   \beta(x)  \cdot  \g(p_* q^* (y)) $.
Hence $h_y \circ S_* \sim h_{S^* y}$.
\>{lem}
 
 \prf  The ``hence'' follows immediately from the main statement, using the
 definition of $h_y$ and the fact that $q_*$ preserves degrees, while $p^*$ multiplies them by $d$.   To begin with:
 
 {\bf Claim.} Let $D' \in NS(V')$, $D = q^* D'$.  Then  $p^* p_* D = \sum_{h \in H} h ^* D$

\prf  Let $Ram_p$ be the ramification divisor of $p: S \to V$, viewed as a divisor on $S$.
Then $q(Ram_p)$ is a proper subvariety of $V'$, and we may choose a representative of the linear
equivalence class of $D'$ whose support has no component contained in $q(Ram_p)$.  
We show equality of the Cartier divisors $p^* p_* D$ and $\sum_{h \in H} h ^* D$.  Since
$S$ is normal, equality as Weil divisors suffices.  
Since $p$ is finite, neither of these has a component of the support contained in
$Ram_p$; and to show equality we may work away from $Ram_p$.  There, $p$ is \'etale,
and $S$ is smooth.  Since $p$ is Galois, the equality is clear.   \eprf

 We explain below, in sequence, the following chain of equalities and $\sim$: 
 $$ \beta(q_* p^* x) \cdot \g (y) = \beta(p^* x) \cdot q^* \g (y) \sim
 {d \inv}  \beta(p^*x) \cdot \sum_{h \in H} h q^* \g(y) = $$
$$=  {d \inv}   \beta(p^*x) \cdot p^* p_* q^* \g (y) = 
 {d \inv}    p_* \beta(p^*x) \cdot p_* q^* \g(y) = $$
$$= {d \inv}   \beta (p_* p^*x) \cdot p_* q^* \g(y) \sim   \beta(x)  \cdot  \g(p_* q^* (y)) $$

\<{itemize}
\item By \eqref{H3} applied to $\xi=p^*x$. 
\item Here we use the fact that $H$ acts by automorphisms on $\bS$, and so respects
the intersection product; thus $u \cdot v = d \inv \sum_{h \in H} (hu) \cdot (hv)$;
and $h \beta p^* x = \beta h p^* x = \beta p^* x$.  
\item By the claim.
\item By the projection formula \eqref{projection}.
\item  By \eqref{H2.1}.
\item While $\g$ is not defined on $NS(S)$, $\rho$ is defined on $V,S,S'$ and commutes
with $p_*$ and with $q^*$.   Hence $\rho$ commutes with $p_*q^*$, so 
$\rho(p_*q^* \g(y)) = p_*q^* (y)$, and hence by Lemma \ref{height1} we have
$\bt(x) \cdot (p_*q^* \g(y)) \sim \bt(x) \cdot \gamma( p_*q^* (y))$.  We aso use
$p_*p^*  x =d x$ on $C_0(V)$.
 \>{itemize}
 
 \eprf

For the sake of a later observation (\ref{alg}), we  record an easier case:

\<{lem}  \vlabel{transpose0}  Assume $p$ is finite above a Zariski open neigborhood
of the support of $\bba \in C_1(\bV)$.  Then
  $$h_{\bD}  \circ S_* (a)=   h_{\bS^* ( \bD)}(a) $$
 \>{lem}
 
 \prf       In this case we have $\bt p^* (a)=p^* \bt(a)$,  since  using the 
  finiteness assumption on $\bS \to \bV$ near $a$, there will be no components of $p^* \bt(a)$ projecting to a point of $C$.    The projection formula then immediately gives the lemma. \eprf
 
%In $A_1(\bV)$ we have $[S_*(a)]=  \rho( q_*(p^*\bba \cdot \bS)) =  q_*(p^*\bba \cdot \bS)$,
%  since  using the 
%  finiteness assumption on $\bS \to \bV$ near $a$, there will be no components projecting to a point of $C$.    
%  Hence 
%  $$ deg(S_*(a)) h_D ( S_*(a)) = \  q_*(p^*\bba \cdot \bS)   \cdot \bD$$
% 
%                                                                                      
% In this last  expression, the first $\cdot$ refers to intersection in  $\bV \times_C \bV$.
% the second to intersection in $\bV$.  By the projection formula, 
% $$q_*(p^*\bba \cdot \bS) \cdot \bD =  p^*\bba \cdot \bS \cdot q^* \bD  $$

% Applying $p_*$ (which preserves degrees of 0-cycles),  and    using the projection
% formula again, we obtain:
%  $$q_*(p^*\bba \cdot \bS) \cdot \bD = \bba \cdot p_*(\bS \cdot q^* \bD  ) $$
%or
% $$\bS_*(\bba) \cdot \bD = \bba \cdot \bS^t(\bD)$$
% So $h_D(S_*(a)) \deg(S_*(a)) = h_{S^t(D)} (a) \deg(a)$.
% Since $\deg(S_*(a)) = d \deg(a)$, the statement follows. \eprf

\smallskip
\para{\bf \noindent Canonical heights}

Let $E=  NS_\Rr(V) = \Rr \tensor NS(V) $.  Then $S^*$ induces a linear endomorphism $S^* |E $ of $E$.      
We can (uniquely) express $E$ as a direct sum of two $S^*|E$-invariant subspaces $E_-,E_+$, 
such that   every  complex eigenvalue
of $S^*$ on $E_-$ (respectively $E_+$) has   absolute value $\leq 1$ (respectively $>1$).  In particular, $ S^* | E_+$ is invertible; let $s$ denote the inverse.  
 When $d=1$, the eigenvalues on $E_-$ are $0$ and roots of $1$, cf. Lemma \ref{global}. 

For any $e \in E_+$, the sequence $s^n(e) $ approaches $0$ exponentially fast.

Let $\kappa \inv$ be spectral radius of $s$, i.e. 
$\kappa = \min \{ |\a|: \a \in spec(S^*|E), |\a| >1\}$,
where $spec( -)$ is the set of eigenvalues.  So $\kappa > 1$.  

Let $\Lambda = E_+^*$ be the dual space to $E_+$.  Let $\lambda$ be the 
dual linear transformation to $S^*$, i.e. $\lambda(F)=F \circ S^*$.

\<{prop}   \lbl{heights1}
 There
 exists a unique function $h : V(K^{alg}) \to  \Lambda $ (extending to a linear $h: C_0(V) \to \Lambda$) such that:
 
 (1)  For any $e \in E$, the function:
  $a \mapsto h(a)(e)$  represents $h_e$ in  ${\mathcal F}(V)$.   

(2)    For all $a \in  V(K^{alg})$:   $h(S_*(a)) =    \lambda (h (a))$ 

When  $E_+$ contains an ample divisor, $h$ is ``proper" in the sense that the inverse image
of a bounded subset of $E^*$ is a limited subset of $V(K^{alg})$.
 \>{prop}

\prf 

Recall the lemma behind Tate's  canonical heights construction
(cf. e.g. \cite{silverman} Theorem 3.20).   let $X$ be a set, $T: X \to X$ a function, $h: X \to \Rr$, $\kappa \in \Rr, \kappa >1$.  
Assume $h\circ T - \kappa h$ is a bounded function on $X$.  Then there exists a unique
$h': X \to \Rr$ with $h-h'$ bounded, and such that $h' \circ T = \kappa h'$.  (Namely $h' = \lim \kappa^{-n} h \circ T^n $.)   

We will apply this to an appropriate function on
$X=V(K^{alg}) \times Y$ for various compact subsets $Y$ of $E_+$.  
  We write $V$ for $V(K^{alg})$ for the rest of this proof.

 Let  $Bdd(X)$ be the space of functions $\phi: V \times E_+ \to \Rr$,
such that for any compact $Y \subseteq E_+$, $\phi | (V \times Y)$ is a bounded function.  
 
We have a surjective homomorphism $\rho: NS_\Rr(\bV) \to NS_\Rr(V)$. 
Choose a linear map $\gamma: NS_\Rr(V) \to NS_\Rr(\bV)$, with $\rho \circ \gamma = Id_V$.  
 We use
the fact that $NS(V)$ is a finite dimensional space; $E=  NS_\Rr(V)$ has a uniquely determined topology of a real topological vector space.

Define $h_0: C_0(V) \times NS(V) \to \Rr$ by $h_0(a,e) = h_{\gamma(e)}(a)$, and extend to $C_0(V) \times NS_\Rr(V)$ by linearity in the second variable.    For any
fixed $e \in NS(V)$, $\rho \bS^t(\gamma(e)) = \rho \gamma(S^t(e))$.  Hence
by Lemma \ref{height1}, $h_{ \bS^t(\gamma(e))} - h_{\gamma(S^t(e))}$ is bounded on $C_0(V)$.
Now $h_0(S_*(a),e)) - d \inv h_{\bS^t(\bD)}(a)$ is bounded as a function of $a$, while
by definition $h_0(a,S^t(D)) = h_{\gamma(S^t(e))} (a)$.  Hence for fixed $e$,
$\delta(a,e) = h_0(S_*(a),e) - d \inv h_0(a,S^t(e))$ is bounded on $C_0(V)$.  Let
$e_1,\ldots,e_r \in NS(V)$ be a basis for $NS_\Rr(V)$.  If $Y \subset E$ is compact  then for some $B \in \Rr$, any $y \in Y$ can be written $y=\sum \a_i e_i$ with $|\a_i| \leq B$;
so $\delta(a,y) = \sum_i \a_i \delta(a,e_i)$ is bounded on $C_0(V) \times Y$.  

In particular, restricting to $X =C_0(V) \times E_+$, we see that 
$h_0(S_*(a),e) - d \inv h_0(a,S^t(e))$ lies in $Bdd(X)$.

Choose $\kappa_1$ with $1 < \kappa_1 < \kappa$.   Define
$T: X \to X$ by $T(a,e) = (S_*(a), \kappa _1 s(e))$.  The spectral radius of
$\kappa_1 s(e)$ is $\kappa_1 \kappa \inv < 1$, so that any orbit of $\kappa_1 s$ is bounded
(since it approaches $0$), and more generally any compact subset of $E_+$ is contained in a
compact, $\kappa_1 s$-invariant subset.
  On the other hand modulo $Bdd(X)$ we have:
$$h_0(T(a,e)) = \kappa_1 h_0( S_*(a),  s(e)) =  \kappa_1 d \inv h_0(a,S^t(s(e))) = \kappa_1   h_0(a,e)$$
  Since $\kappa_1>1$, the Tate lemma applies on $C_0(V) \times Y$ for any
  compact, $\kappa_1 s$-invariant subset $Y$ of $E_+$.  
 So  for each such $Y$    there exists a unique $h_1^Y: C_0(V) \times Y \to \Rr$ with
 $h_1^Y-h_0$ bounded and $h_1^Y(T(a,e))=\kappa_1 h_1^Y(a,e)$.   Hence
   there exists a unique $h_1: X \to \Rr$
with $h_1 - h_0 \in Bdd(X)$, and $h_1(T(a,e))=\kappa_1 h_1(a,e)$.    For any fixed $a$, the function $h_1(a,e)$ is defined as a limit
of linear functions in $e$, so $h_1(a,e)$ is linear in $e$.  We have
$h_1(S_*(a), e ) = h_1( T(a,\kappa_1 \inv   S^* (e))) = \kappa_1 h_1(a,\kappa_1 \inv S^*(e)) = h_1(a,S^*(e))$.    
  
  Let $h(a)$ be the linear map:  $e \mapsto h_1(a,e)$.  Then (1),(2) are clear.  The remark on the ample divisor is also clear.

\eprf

 Consider now the case   $d=1$.  When is
    $E_+ \neq 0$?  
Recall $d' = \deg(q)$, the forward degree of $S$.   
\<{lem}  If $d'>1$ then $E_+ \neq 0$.  \>{lem}

\prf We have on $E$ an $n$-multilinear form, the intersection product.  Write
$e^n$ for $e \cdot \cdots \cdot e$.   If $e \in E$ is the class
of an ample divisor, then $e^n >0$, and using the projection formula
$$((S^*)^m e )^n = (d')^{m} e^n$$
As this grows exponentially in $m$, 
the lemma   follows from the following fact from linear algebra:

\claim{1}  Let $V$ be a finite dimensional complex vector space, and let $m: V^k \to \Cc$
be a multilinear map.  Let $v \in V$.  Let $T \in End(V)$, and suppose every eigenvalue of 
$T$ has absolute value $\leq 1$.  Then $|m(T^av,\ldots,T^av)|$ is bounded
by a polynomial in $a$.   
\prf  Let $v_1,\ldots,v_n$ be a basis with respect to which $T$ has Jordan normal form.
Let $v= \sum \a_i v_i$.  Then $T^a v = \sum q_i(a) v_i$ where $|q_i(a)|$ is polynomially
bounded in $a$.  It follows from multilinearity that $m(T^av,\ldots,T^av)$  is polynomially bounded too. \eprf

  \eprf

We thus obtain a nontrivial function into a finite dimensional real dynamical system whenever
the algebraic dynamics is non-birational.    But even in the birational case, it would seem to be rare that $E_+$ is trivial, and it would be nice to obtain more geometric information.  We have
the following from global linear algebra, showing that $E_+$ is trivial only when 
$S^*$ is essentially uni-by-nilpotent on $E$.  Here we use the fact that the action of $S^*$ on $E_+$ arises from the action on a finitely generated group, the N\'eron-Severi group $N$ of $V$; the lemma would not be true if $N$ were allowed to be an arbitrary $\Qq$-space.

\<{lem} \lbl{global}    Let $N$ be a finitely generated Abelian group, and $T \in End(N)$.   Let $V=N \tensor \Cc \cong \Cc^n$, and write $T$ for $T \tensor \Cc$.
Assume
  every   eigenvalue of $T$ has
   complex absolute value $\leq 1$.
Then every such eigenvalue is $0$ or a root of $1$.  \>{lem}

\prf    These eigenvalues
lie in some number field $L$, and form a set closed under conjugation. For each non-archimedean absolute value $p$ of $L$, each eigenvalue has absolute value $\leq 1$,
since the compact open set $L_p \tensor N$ is left invariant by $T$.  Thus each eigenvalue is an algebraic integer, all of whose conjugates have complex norm $\leq 1$;
by the product formula, it is $0$ or else every conjugate has complex norm $1$; in the
latter case it is a root of unity.
 \eprf

 Hence if $d=1$; then $E_+$ is the image of $(S^t)^n ((S^t)^n-1)^m$, for large enough $m,n$. 
 
{\em Proof of Theorem \ref{thmh}} (a) has been proved above.  As for (b), by assumption, and by Lemma \ref{global}, every eigenvalue of $\phi^*$ is either $0$ or of absolute
value $>1$.   Let $D$ be an ample divisor.  Then we may write $D=D_1+D_2$
where $(S^*)^m D_1 = 0$ for some $m$, and $D_2 \in E_+$.   Assume 
 $\lambda(a)=0$.  Then $\lambda(\phi^k(a))=0$ for all $k$, so    $h_{D_2}(\phi^k(a))$ is bounded uniformly in $k$.   On the other hand $h_{D_1} ( \phi^m(b)) \sim 0$  
 since $h_{D_1} \circ \phi^m \sim h_{(S^*)^m D_1} $ but $ (S^*)^k D_1 = 0$.  This holds
 uniformly for all $b$, in particular for $b =\phi^k(a)$, so  $h_{D_1} ( \phi^{m+k}(a)) \sim 0$  
uniformly in $k$.   Thus $h_D(\phi^k(a)) \sim 0$.  
% Thus for large enough $m$ we have $(S^*)^m E = E_+$.  
  \qed

{\em \bf  \noindent Examples} 

(1)   If $V$ is a curve, then $NS(V)$ is one-dimensional, so $S$ acts by   multiplication by a scalar.  In the case of a curve the scalar is   the degree of $\pi_2 $ on $S$. 
Thus the condition (*) in this case is that the degree of $\pi_2$ is greater than that of $\pi_1$.
The classical case has $d=1$, and the condition is the $\pi_2$ has degree $>1$.   
 
(2)   If $V= \Pp^n$  then $NS(V)$ is one-dimensional, and
 $S$ acts by multiplication by a scalar, which is positive unless $S$ is linear.  

(3)   Suppose $S$ is the graph of a morphism $s$.  Call and Silverman \cite{CS} assume the existence of $e \in E$ with $S^*(e) = \kappa e$, $
\kappa >1$.  In this case we have  $e \in E_+$, and the Call-Silverman canonical height
$h_{V,e,s}$ is given by $h_{V,e,s} (a)  = h(a)(e)$.
 
 More generally if all the  eigenvalues of $S^* | E_+$ are real, then $h$ is captured by 
 $\dim(E_+)$ canonical heights into $\Rr$.  But it is easy to find examples with no
 real eigenvalues.  ( Using $SL_n$ actions on compactifications of $(G_m)^n$.)

\para{\noindent \bf Generic algebraicity}  For simplicity, assume $S$ is the graph of a morphism $\phi$, 
 and that $S^*(D) = \lambda D$ for some ample divisor $D$ and some $\lambda >1$,  so that we have a canonical height function $H_D: V(K^{alg}) \to \Rr$.   For Abelian varieties over number fields,  Silverman suggested  that the canonical height may be transcendental for sufficiently general algebraic points.  We note here that for many (if not all) dynamics over function fields, the canonical height of sufficiently irrational algebraic points is an algebraic number.
 
\def\bE{{\bf E}}
Let $\bE$ be the blowup locus of $\bS \to \bV$.  $\bE$ projects to a finite subset  $\{c_1,\ldots,c_k\}$ of $C$, so $\bE= \union_i E_i$,
with $E_i = E_{c_i}$.  
Even for non-isotrivial dynamics, $\bE$ may be empty.  In the case of dynamics
on curves it is at most finite.   Note: 

\<{lem}  Assume $\phi$ is a rational map $V \to V$, defined 
outside a finite set $E$,  over a field $k$.  Let $l$ be a prime
bigger than $|E|$ and $\deg(\phi)$.   Let $d \in V(k^{alg})  $ with $l | [k(d):k]$.
Then $\phi^m(d)$ is defined for any $m \in \Nn$.
 \>{lem}

\prf  If $e \in E$ then $[k(e):k] \leq |E| < l$; hence $d \notin E$, so $\phi(d)$
is defined.  We have $[k(d,\phi(d)):k(\phi(d))] \leq \deg(\phi)$, so
if $[k(\phi(d)):k]$ is prime to $l$ then so is $[k(d,\phi(d)):k]$ and hence $[k(d):k]$,
a contradiction.  Thus $[k(\phi(d)):k]$ is also divisible by $l$.  We continue inductively.
\eprf

Let $\tV$ be 
the set of algebraic points $a$ of $V$, such that curve $\rho(\phi^n(a))$ 
corresponding to $\phi^n(a)$  does not meet $\bE$, for any $n$.  Let $c_1,
\ldots,c_n$ be the support of the projection of $\bE$ to $C$, $V_i = V_{c_i}$,
$\phi_i$ the restriction of $S$ to $V_i$.  Assume $\bE$ is finite, so that $\phi_i$ is a rational map (defined away from $E_i$.)  Assume the base field $k$ is not 
real closed or algebraically closed.  (The heights can be evaluated with respect to $k^{alg}(C)$;
the points in question will be in $k(C) \m k^{alg}$.)

and let $d_i \in V_i(k^{alg})$ be such that
$[k(d_i,c_i): k(c_i)]$ is divisible by a large prime $l$.  By the lemma, $\phi_i^m(d_i)$
is defined for all $m$.  By the approximation lemma, choose $a \in \bV(K^{alg})$   such that under the map
$\res: \bV(K^{alg}) \to V_(k^{alg})$ associated with the valuation at $c_i$, 
we have $\res(a)=d_i$.  Then $a \in \tV$.  
In this sense, `sufficiently general' algebraic points lie in $\tV$; meaning that
the residue at a finite number of places is to be highly irrational, i.e. of high degree.

Presumably  this is true in general, and not only when $\bE$ is finite.

\<{lem} \vlabel{alg}  Let $a \in \tV$.  Then $H_D(a)$ is an algebraic number.   \>{lem}
 
  \prf   Let $ f=\deg(a)$.  
 Let $\bS,\bD$ be subvarieties of $\bV \times_C \bV,\bV$ defined over $k$
and restricting to $S,D$, as above.
%  Let  $\bD_m =  (\bS^t)^m (\bD)$, and $d_m = \beta(a) \cdot \bD_m$.  
Let $d_m = (\bS_*^m(\bba)) \cdot \bD $. 
 By   Lemma \ref{transpose0} we have $d_m= \bba \cdot (\bS^t)^m (\bD)$, and 
 $h_{\bD} ( \phi^m(a) )= f \inv d_m$.
 %$h_{\bD}(\bS_*^m(a)) =   h_{(\bS^t)^m (\bD)} (a)=   f \inv d_m$. 
So $H_D(a) =   \lim_ {m \to \infty} f \inv \kappa ^{-m}   d_m$.   We will prove
 the stronger claim, that the generating series $\sum d_m t^m$ is rational with algebraic
 coefficients.  
  
   We use the finite-dimensionality of $W=NS(\bV)$.  $\bS^t$ acts on $NS(\bV)$ linearly.
 We have a linear map $l: End(W) \to \Rr$, namely $w \mapsto \beta(a) \cdot w$.  It
 follows by using Jordan form over $\Qq^{alg}$
   that 
$d_m = \sum p_i(m) \kappa_i^m$ for some polynomials $p_i$ with algebraic coefficients,
and some $\kappa_i \in \Qq^{alg}$ and $p_i \in \Qq^{alg}[T]$.  It 
 follows that $\lim_{m \to \infty} f \inv d_m \kappa^{-m}$
(being finite) is algebraic.  \smallskip

If for example $V$ is the generic point of an Abelian scheme  over $C$, with dynamics given by multiplication,   there will be no blowup locus, and the canonical height will be algebraic everywhere.

\bigskip \noindent
Current addresses:

\bigskip\noindent
UFR de Math\'ematiques \par\noindent
Universit\'e Paris 7 - Case 7012
\par\noindent 2, place Jussieu
\par\noindent
75251 Paris Cedex 05 \par\noindent
France \par\noindent
e-mail: {\tt zoe@logique.jussieu.fr}

\bigskip\noindent
Dept. of Mathematics, Hebrew University, Jerusalem, Israel; \par\noindent
(2008) Yale University Mathematics Dept., PO Box 208283 \par\noindent
New Haven, CT 06520-8283, USA\par\noindent
e-mail: {\tt ehud@math.huji.ac.il}

 .
\end{document}